\documentclass[12pt]{amsart}
\usepackage{latexsym}
\usepackage{amssymb}
\usepackage{amsmath}

\newtheorem{thm}{Theorem}[section]
\newtheorem{lem}[thm]{Lemma}
\newtheorem{claim}[thm]{Claim}
\newtheorem{cor}[thm]{Corollary}
\newtheorem{prop}[thm]{Proposition}

\newtheorem{prob}[thm]{Problem}

\theoremstyle{remark}

\numberwithin{equation}{section}

\newcommand{\beq}{\begin{equation}}
\newcommand{\eeq}{\end{equation}}

\def\real{\hbox{\rm\setbox1=\hbox{I}\copy1\kern-.45\wd1 R}}
\def\prob{\hbox{\rm\setbox1=\hbox{I}\copy1\kern-.45\wd1 P}}
\def\natural{\hbox{\rm\setbox1=\hbox{I}\copy1\kern-.45\wd1 N}}
\def\integers{\hbox{\rm\setbox1=\hbox{-}\copy1\kern-1.2\wd1 Z}}

\DeclareMathOperator{\sgn}{sgn}

\newcommand{\B}{\mbox{$\mathcal B$}}

\allowdisplaybreaks

\begin{document}
\title[Coboundary existence]{Existence and Non-existence of Solutions to the Coboundary Equation for Measure Preserving Systems (DRAFT)}
\author{Terry Adams and Joseph Rosenblatt}
\date{June 2019}

\begin{abstract} 
A fundamental question in the field of cohomology of dynamical systems is to determine when there are solutions to the coboundary equation: 
\[
f = g - g \circ T  .
\]
In many cases, $T$ is given to be an ergodic invertible measure preserving transformation 
on a standard probability space $(X, \B, \mu )$ and 
$f: X\to \real$ is contained in $L^p$ for $p \geq 0$.  We extend previous results 
by showing for any measurable $f$ that is non-zero on a set of positive measure, 
the class of measure preserving $T$ with a measurable solution $g$ is meager 
(including the case where $\int_X f d\mu = 0$).

From this fact, a natural question arises: given $f$, does there always exist a solution pair $T$ and $g$? 
In regards to this question, our main results are:
\begin{itemize}
\item Given measurable $f$, 
there exists an ergodic invertible measure preserving transformation $T$ 
and measurable function $g$ such that $f(x) = g(x) - g(Tx)$ for a.e. $x\in X$, if and only if 
$\int_{f > 0} f d\mu = - \int_{f < 0} f d\mu$ (whether finite or $\infty$). 
\item Given mean-zero $f \in L^p$ for $p \geq 1$, 
there exists an ergodic invertible measure preserving 
$T$ and $g \in L^{p-1}$ such that  $f(x) = g(x) - g( Tx )$ for a.e. $x \in X$. 
\item In some sense, the previous existence result is the best possible. 
For $p \geq 1$, there exists a dense $G_{\delta}$ set of mean-zero $f \in L^p$ such that 
for any ergodic invertible measure 
preserving $T$ and any measurable $g$ such that $f(x) = g(x) - g(Tx)$ a.e., 
then $g \notin L^q$ for $q > p - 1$.  
\end{itemize}

Finally, it is shown that we cannot expect finite moments for solutions $g$, when $f \in L^1$. 
In particular, given any $\phi : \real \to \real$ such that $\lim_{x\to \infty} \phi (x) = \infty$, 
there exist mean-zero $f \in L^1$ such that for any solutions $T$ and $g$, 
the transfer function $g$ satisfies: 
\[
\int_{X} \phi \big( | g(x) | \big) d\mu = \infty . 
\]
\end{abstract}
\maketitle

\section{Introduction}\label{intro}

We give new fundamental results concerning solutions to the coboundary equation:
\begin{eqnarray} \label{eqn:cob-1}
f = g - g \circ T . 
\end{eqnarray}
There has been substantial progress in many cases such as homogeneous spaces, 
smooth actions, lie groups, as well as many other important families of dynamical systems. 
Most previous research focuses on the case where a measurable transformation, 
or topological dynamical system is specified, and a solution $g$ is sought 
for individual $f$ or families of functions $f$ (e.g., H\"{o}lder $f$). 
In this paper, we study the situation from the general perspective 
of solutions $T$ and $g$ where $f$ may be any real-valued measurable function, 
or function $f \in L^p$ for $p \geq 0$. 

Let $(X, \B, \mu)$ be a standard probability space
\footnote{In this paper, standard probability space means isomorphic to $[0,1]$ with Lebesgue measure.}, 
and for $p > 0$, define the standard $L^p$ space, 
$L^p = \{ f: X \to \real | f\mbox{ is measurable and } \int_X | f |^p d\mu < \infty \}$. 
For $p\geq 1$, define $L^p_0 = \{ f\in L^p : \int_X f d\mu = 0 \}$. 
Also, $L^{\infty}$ is the set of essentially bounded measurable functions on $(X, \B, \mu)$ 
and similarly, $L^{\infty}_0$ are functions in $L^{\infty}$ with zero integral. 
The space $L^0$ is the set of measurable functions on $(X, \B, \mu)$. 
Let $\mathcal{M}$ be the family of invertible measure preserving transformations 
defined on $(X, \B, \mu)$ and $\mathcal{E}$ is the family 
of ergodic invertible measure preserving transformations on $(X, \B, \mu)$. 
We obtain the following main positive result:
\begin{thm}[Existence of solutions]
\label{pos-thm}
Let $1 \leq p \leq \infty$ and suppose $f \in L^p_0$.  
There exist $T \in \mathcal{E}$ and $g\in L^{p-1}$ 
such that  $f(x) = g(x) - g(Tx)$ for a.e. $x \in X$. 
\end{thm}
In some sense, Theorem \ref{pos-thm} gives the best possible positive result. 
The following theorem demonstrates a major limitation for solutions to the coboundary equation. 
In particular, typically, there is no solution $g$ in the same integrability class as $f$, 
even when allowing $T$ to range over all of $\mathcal{E}$. 
\begin{thm}[$L^q$ non-existence]
\label{neg-thm}
Given $1 \leq p < \infty$, there exist $f \in L^p_0$ such that for any solution 
$T \in \mathcal{E}$ and measurable $g$ to the coboundary equation 
$f = g - g \circ T$, then $g \notin L^q$ for $q > p - 1$. 
More generally, there exists a dense $G_{\delta}$ set 
$\mathcal{G}_p \subset L^p_0$ such that for any $f \in \mathcal{G}_p$, 
and any solution pair $T, g$ with $T \in \mathcal{E}$, then 
$g \notin L^q$ for $q > p - 1$. 
\end{thm}
The solution $g$ is referred to as the transfer function for coboundary $f$. 
Theorem \ref{neg-thm} implies for generic mean-zero $f \in L^p$ 
for $p < 2$, that any transfer function is not integrable, regardless of $T \in \mathcal{E}$. 
However, for $f \in L^1_0$, we can always find a solution with measurable $g \in L^0$. 
For the case where $f$ is only assumed to be measurable, 
we give a straightforward equivalent condition for the existence of a measurable transfer function.  
Also, Theorem \ref{meas-thm} highlights the need to control $T$, or the inter-dependence 
of $T$ and $f$, if one hopes to find a measurable transfer function. 
\begin{thm}[Measurable transfer functions]
\label{meas-thm}
Suppose $(X, \B, \mu)$ is a standard probability space 
and $f \in L^0$ is non-zero on a set of positive measure. 
\begin{itemize}
\item The class of ergodic invertible measure preserving transformations 
$T$ such that $f = g - g \circ T$ has a measurable solution $g$ is first category (i.e., meager); 
\item The coboundary equation $f = g - g\circ T$ has a solution pair, $T \in \mathcal{E}$, $g \in L^0$, 
if and only if $\int_{f > 0} f d\mu = - \int_{f < 0} f d\mu$, whether both 
integrals are $\infty$ or finite.  This is an extension of Anosov's observation 
\cite{Anosov73} to include the case where $f$ is not integrable. 
\end{itemize}
\end{thm}

\section{Connections to Previous Research}
There has been substantial interest in the study of the cohomology of dynamical systems.  
Much of the recent focus is on smooth dynamics including hyperbolic actions 
or actions of lie groups.  Powerful rigidity or local rigidity results have been 
obtained involving cocycles. Some of the earliest results include \cite{KS97} and \cite{KS94}. 
Cocycle rigidity depends closely on solving the coboundary equation, 
since the difference between cohomologous cocycles is a coboundary. 
Liv\u{s}ic \cite{Liv72} provided one of the earliest regularity results in this setting 
by demonstrating H\"{o}lder cocycle rigidity for families of U-systems, 
topological Markov chains and Smale systems. 
More recently, this H\"{o}lder regularity has been extended 
to nonuniformly expanding Markov maps \cite{GO2005}, and 
to Weyl chamber flows or twisted Weyl chamber flows \cite{Vin2017}. 

We will consider the coboundary equation in a general context. 
In the setting of topological dynamics, 
the following was observed in Gottschalk and Hedlund \cite{GH55} 
and later extended by Browder \cite{Brow58}:
a bounded continuous function $f$ is a coboundary 
for a homeomorphism on a compact space 
if and only if the following is uniformly bounded for positive $n$,
\[
| \sum_{i=0}^{n-1} f ( T^i x ) | . 
\]

\subsection{Schmidt's condition}
\label{schmidt}
The following associated condition for measurable dynamics can be found 
in \cite{Hal76} and \cite{Schmidt77}. 
A measurable function $f$ is a coboundary for $T\in \mathcal{E}$ 
if and only if for each $\delta > 0$, 
there exists $M_{\delta} \in \natural$ such that for $n \in \natural$, 
\[
\mu \big( \{ x \in X : | \sum_{i=0}^{n-1} f (T^i x ) | \leq M_{\delta} \} \big) > 1 - \delta . 
\]
This condition will be used in section \ref{cat-sec} to show for any 
measurable function $f$ that is essentially non-zero, 
then the class of ergodic invertible measure preserving transformations 
$T$ such that  $f = g - g\circ T$ has a measurable solution $g$ is meager 
(first category). 
Anosov \cite[Theorem 1]{Anosov73} 
demonstrated that there are no measurable solutions $g$ 
in the case that $f$ is integrable and $\int_{X} f d\mu \neq 0$. 
However, our category results apply in the situation that $\int_{X} f d\mu = 0$. 

\subsection{Non-measurable solutions}
Using the axiom of choice, we can always obtain a solution $g$. 
Partition $X$ into orbits.  For each orbit $\mathcal{O}$, 
choose a single point $x_0 \in \mathcal{O}$. 
The coboundary equation leads to the following telescoping series, for $n > 0$, 
\[
g( T^n x ) = g(x) - \sum_{i=0}^{n-1} f ( T^i x ) , 
\]
and for backward iterates, 
\[
g( T^{-n} x ) = g(x) + \sum_{i=1}^{n} f ( T^{-i} x ) . 
\]
If we define $g(x_0) = 0$, then the recursion formulas above 
uniquely determine $g$ at all points along the orbit and at a.e. $x\in X$. 
However, the result of Anosov implies this $g$ is not measurable 
when $f$ has a non-zero integral. 

Here is another case where this construction clearly leads to a non-measurable solution. 
Suppose $\alpha$ is irrational and $0 < \alpha < 1$. 
Define $f$ on $[0,1]$ by:
\begin{eqnarray*} 
f(x)= 
\left\{\begin{array}{ll}
\alpha, & \mbox{if $x \leq \frac{1}{1+\alpha}$}, \\ 
- 1 & \mbox{if $x > \frac{1}{1+\alpha}$}.
\end{array}
\right.
\end{eqnarray*}
The integral of $f$ is zero.  Since $g(x) = 0$ for a single point in each orbit, 
then the space $X$ equals the following disjoint union (modulo measure zero sets), 
\[
\bigcup_{i = -\infty}^{\infty} T^i \big( \{ x\in X : g(x) = 0 \} \big) . 
\]
Since $T$ is measure preserving, the set $\{ x\in X : g(x) = 0 \}$ is not measurable 
and consequently, $g$ is not measurable. 

There are other cases where it is known that the coboundary equation has 
no measurable solution $g$.  It was pointed out in \cite{Hal76} 
that if $f$ is a non-trivial step function taking on two values, then the transformation $T$ 
must have a non-trivial eigenvalue. 
Thus, if $T$ is weakly mixing and $f$ is a 2-step function, 
there is no measurable solution $g$. 

\subsection{Bounded coboundaries}
This raises the question of when do solutions exist for classes of measurable 
functions $f$, when $T$ is allowed to range over $\mathcal{E}$. 
In \cite{TAJR}, it is shown that any finite step, mean-zero function is a coboundary 
for some ergodic invertible measure preserving transformation with a bounded 
transfer function $g$. 
In particular, $T$ may be chosen in one of the following categories: 
\begin{enumerate}
\item $T$ is a transformation with discrete spectrum;
\item $T$ is a product of rotations;
\item $T$ is a finite extension of a product of rotations.
\end{enumerate}
Also, in \cite{TAJR}, the existence of solutions is extended to mean-zero bounded functions. 
The case of general $L^p_0$ functions is more subtle and addressed in this paper. 

The paper \cite{Kw} partially addresses the case of bounded coboundaries.  However, the arguments 
given in \cite{Kw} are viewed as containing a gap, 
and the main theorem does not apply in general beyond 
the case of continuous functions $f$.  

\subsection{Operator viewpoint}
The coboundary equation has been viewed from the perspective 
of operator theory.  Note that the coboundary equation may be written as,
\[
f = (I - U_T) g 
\]
where $U_T$ is the Koopman operator defined by $U_T(g) = g\circ T$, 
and $I$ represents the identity operator. 
Study of the operator $( I - T )$ when $T$ is a linear operator 
(and not necessarily unitary) 
goes back to the $19^{th}$ century \cite{Neumann1877}. 
Similar to the case of real or complex numbers, 
for an operator $T$ with norm $| T | < 1$, 
then $I - T$ has an inverse and 
\[
\big( I - T \big)^{-1} = \sum_{i=0}^{\infty} T^k . 
\]
However, for measure preserving transformations, 
$| U_T | = 1$, and solving $f = (I - U_T) g$ becomes more complicated. 
Iterative techniques were given in \cite{Dot69, Dot70, Dot71, Gro76} 
as an aid for solving the coboundary equation in this setting. 
The paper \cite{LS} shows that for a given $T$, when a solution 
exists, it may be obtained in closed form as the following point-wise limit a.e.:
\[
g(x) = \lim_{n\to \infty} \frac{1}{n} \sum_{k=1}^{n} \sum_{i=0}^{n-1} f ( T^i x ) . 
\]
Also, the authors extend their results from the classical Poisson equation, 
$f = (I - U_T) g$ to the case of fractional coboundaries \cite{DL2001}. 
Their main results produce equivalent conditions for solutions 
to occur for fixed $T$.  

Our main results can be recast in terms of operators in the following way. 
\begin{cor}[Operator theoretic statement of Theorems \ref{pos-thm} and \ref{neg-thm}]
Let $(X, \B, \mu)$ be a standard probability space and 
$\mathcal{E}$ be the set of all ergodic invertible measure preserving 
transformations on $(X, \B, \mu)$. 
Then Theorems \ref{pos-thm} and \ref{neg-thm} are equivalent 
to the following statements respectively, 
\[
L^p_0 \subset \bigcup_{T\in \mathcal{E}} \Big( I - U_T \Big) \Big( L^{p-1} \Big) 
\]
and 
\[
L^p_0 \cap \bigcup_{T\in \mathcal{E}} \bigcup_{q > p - 1} \Big( I - U_T \Big) \Big( L^{q} \Big)
\ \mbox{is meager in}\ \ L^p_0 . 
\]
\end{cor}

\subsection{Ergodic averages}
One of the main applications of coboundary solutions is to find functions 
for which the ergodic averages are controlled and converge rapidly. 
In the case where $f$ is a coboundary for $T$ with integrable transfer function 
$g$, then all moving averages $(v_n, L_n)$ converge pointwise 
for an increasing sequence $L_n \in \natural$, 
\[
\frac{1}{L_n} \sum_{i=1}^{L_n} f ( T^{v_n + i} x ) \to \int_X f d\mu . 
\]
Other results \cite{VW2004} characterize the rate of convergence of $L^{\infty}$ functions 
using approximation by coboundaries where the transfer function lands 
in a specific $L^p$ space.  For $p\geq 1$, the rate is on the order 
of $n^{-p}$. 
For stationary processes exhibiting randomness 
(e.g., positive entropy, random fields), there is 
a technique for decomposing the process into coboundary and martingale 
components.  See \cite{G69,Vol93, G2009, KKM2017, G2018} and 
the references contained therein for background on this technique 
and its applications. This has made it possible to establish common 
statistical laws (central limit theorem, weak invariance principle) 
in these cases. 

\subsection{Nonsingular transformations}
There is also extensive research on the connections of coboundaries 
to nonsingular transformations.  We do not discuss this in detail, 
but encourage the interested reader to check 
\cite{Aar97, DS2012} for its connections, including 
the existence of equivalent finite or sigma-finite invariant measures. 

\section{\bf Coboundary Existence Theorem}\label{COB}
In this section, we prove Theorem \ref{pos-thm}, 
although it is restated here in an equivalent form. 
We will also show later that this is generally the best possible result.
\medskip

\begin{thm}
\label{cet}
Let $p\in \real$ be such that $p \geq 1$.  Given any mean zero function $f\in L^p$, there exists
an ergodic measure preserving dynamical system $(X,\B,\mu, T)$ and a function $g \in L^{p-1}$
such that $f(x) = g(Tx) - g(x)$ for almost every $x\in X$. 
\end{thm}
For the case of $L^\infty$, this theorem follows from the results of \cite{TAJR}.
However, \cite{TAJR} did not handle unbounded functions.
The technique given here is more straight-forward and can be adapted
to find ergodic measure preserving transformations for unbounded functions. 

First, we define balanced partitions and balanced uniform towers, as was defined in \cite{TAJR}.
Then we state and prove lemmas modified from ones given in \cite{TAJR}.
These are used in a new construction to establish Theorem \ref{cet}.

\subsection{Coboundary Extensions}
In this section, we show how to extend a coboundary for an induced transformation 
to a coboundary for the full transformation.  
Let $T: X \to X$ be an ergodic measure preserving transformation.  
Let $A \subset X$ be a set of positive measure.  
Suppose 
\[
T_A (x) = T^{n_A(x)} (x) , x\in A, 
\]
is the induced transformation defined on $A$.  
See \cite{Pet83} for the definition of an induced transformation. 
Given measurable function $f: X \to \real$ and $x \in A$, define 
\[
f_A (x) = \sum_{i=0}^{n_A(x) - 1} f ( T^i x ) . 
\]
We have the following lemma which will be used to prove Theorem \ref{cet}. 
\begin{lem}
\label{cob-ext-lem}
Let $f: X \to \real$ be a measurable function.  Suppose $f_A$ is a coboundary 
for induced transformation $T_A$ with transfer function $g_A$ such that 
$f_A = g_A \circ T - g_A$.  Then $f$ is a coboundary for transformation $T$ with 
transfer function $g$ defined such that for $x \in A$ and $0\leq j < n_A(x)$, 
\[
g( T^j x ) = g_A ( x ) + \sum_{i=0}^{j-1} f ( T^i x ) . 
\]
In particular, 
\[
f(x) = g ( T x ) - g ( x ) . 
\]
\end{lem}
\begin{pf}
Let $x \in A$.  First, suppose $0\leq j < n_A(x) - 1$. 
Thus, 
\[
g ( T T^j x ) - g ( T^j x ) = g_A(x) + \sum_{i=0}^{j} f ( T^i x ) - g_A(x) - \sum_{i=0}^{j-1} f ( T^i x ) 
= f ( T^j x ) . 
\]
If $y = T^j x$, then $f(y) = g ( Ty ) - g ( y )$. 
Now suppose $y = T^{n_A(x) - 1} x$. 
Then 
\begin{eqnarray}
g ( T y ) - g ( y ) &=& g_A ( T_A x ) - \big( g_A ( x ) + \sum_{i=0}^{n_A(x)-2} f ( T^i x ) \big) \\ 
&=& f_A ( x ) - \sum_{i=0}^{n_A(x)-2} f ( T^i x ) \\ 
&=& \sum_{i=0}^{n_A(x)-1} f ( T^i x ) - \sum_{i=0}^{n_A(x)-2} f ( T^i x ) \\
&=& f ( T^{n_A(x) - 1} x ) = f ( y ) . 
\end{eqnarray}
This proves that $f$ is a coboundary for $T$ with transfer function $g$ for almost every $y \in X$. 
$\Box$
\end{pf}

\subsection{Tower Constructions for Finite-step and Bounded Functions}
\begin{lem}
\label{two-step-lem}
Suppose $A \subset X$ has positive measure and $f: A \to \real$ is contained 
in $L^{\infty}_0$ and takes on 2 steps. 
Given $h \in \natural$ and $\epsilon > 0$, there exist $h_1, h_2 > h$, 
disjoint $I_1, I_2 \subseteq A$ and an invertible measure preserving map $T$ such that: 
\begin{eqnarray}
\mu \big( \bigcup_{i=0}^{h_1 -1} T^i I_1 \cup \bigcup_{i=0}^{h_2 -1} T^i I_2 \big) & = & \mu (A) , \\
T^i I_1, 0\leq i < h_1, && T^i I_2, 0\leq i < h_2\ \mbox{ are all disjoint}, \\
|\sum\limits_{i=0}^{k} f(T^i x) | & \leq & \| f \|_{\infty} \ \ \mbox{ for } x\in I_j, k < h_j , j=1,2, \label{nbig4}\\
| \sum\limits_{i=0}^{h_j -1} f(T^i x) | & < & \epsilon \ \ \mbox{ for } x\in I_j , j=1,2, \label{nbig5} \\
\sum\limits_{i=0}^{h_j -1} f(T^i x) &=&  \sum\limits_{i=0}^{h_j -1} f(T^i y) \ \ \mbox{ for } x, y\in I_j , j=1,2, \label{nbig6} \\ 
1 - \epsilon < \frac{h_1}{h_2} & < & 1 + \epsilon . 
\end{eqnarray}
\end{lem}

\begin{pf}
WLOG, assume $A = [0,1]$. 
Suppose $f = b I_B - c I_C$ is mean zero for $b, c > 0$ and disjoint $B, C$ such that $B \cup C = A$. 
The case where ${b}/{c}$ is rational is straightforward, so we assume 
${b} / {c}$ is irrational. 
There exist $\delta_1, \delta_2$ of the same sign, and $p_1, q_1, p_2, q_2$ 
such that $| \delta_2 | < | \delta_1 | < \epsilon$, 
$p_1 < \epsilon p_2$, $q_1 < \epsilon q_2$, 
$p_2 + q_2 - p_1 - q_1 > h_1$, $p_1 b - q_1 c = \delta_1$ and 
$p_2 b - q_2 c = \delta_2$. 
WLOG, assume $0 < \delta_2 < \delta_1 < \epsilon$.  The case 
where $\delta_1, \delta_2$ are negative follows similarly. 
Let $p_3 = p_2 - p_1$ and $q_3 = q_2 - q_1$. 
Note, 
\[
p_3 b - q_3 c = \delta_2 - \delta_1 < 0 . 
\]
Let $\delta_3 = \delta_1 - \delta_2$. 
Split $B$ into two disjoint sets $B_1, B_2$ such that 
\begin{eqnarray}
\mu (B_1) = \frac{ p_2 \delta_3 }{ ( p_2 + q_2 )\delta_3 + ( p_3 + q_3 )\delta_2 } 
&\mbox{and}& 
\mu (B_2) = \frac{ p_3 \delta_2 }{ ( p_2 + q_2 )\delta_3 + ( p_3 + q_3 )\delta_2 } . 
\end{eqnarray}
Note, 
\begin{eqnarray}
\mu (B_1) + \mu (B_2) &=& 
\frac{ p_2 \delta_3 + p_3 \delta_2 }{ ( p_2 + q_2 )\delta_3 + ( p_3 + q_3 )\delta_2 } \\
&=& \frac{ p_2 (q_3 c - p_3 b) + p_3 (p_2 b - q_2 c) }{ ( p_2 + q_2 )( q_3 c - p_3 b) + ( p_3 + q_3 )(p_2 b - q_2 c)} \\ 
&=& \frac{ (p_2 q_3 - p_3 q_2) c }{ ( p_2 q_3 - p_3 q_2 ) b + (p_2 q_3 - p_3 q_2) c} \\ 
&=& \frac{ c }{ b + c } = \mu (B) . 
\end{eqnarray}
Similarly, split $C = C_1 \cup C_2$ such that 
\begin{eqnarray}
\mu (C_1) = \frac{ q_2 \delta_3 }{ ( p_2 + q_2 )\delta_3 + ( p_3 + q_3 )\delta_2 } 
&\mbox{and}& 
\mu (C_2) = \frac{ q_3 \delta_2 }{ ( p_2 + q_2 )\delta_3 + ( p_3 + q_3 )\delta_2 } . 
\end{eqnarray}
Divide $B_1$ into $p_2$ disjoint sets $B_{1,j}$ for $j \in \{ 1, 2, \ldots , p_2 \}$. 
Divide $C_1$ into $q_2$ disjoint sets $C_{1,j}$ for $j \in \{ 1, 2, \ldots , q_2 \}$. 
Divide $B_2$ into $p_3$ disjoint sets $B_{2,j}$ for $j \in \{ 1, 2, \ldots , p_3 \}$. 
Divide $C_2$ into $q_3$ disjoint sets $C_{2,j}$ for $j \in \{ 1, 2, \ldots , q_3 \}$. 
Thus, $\mu (C_{1, j}) = \mu (B_{1,k})$ for $j \in \{ 1,2, \ldots , p_2 \}$ 
and $k \in \{ 1,2, \ldots , q_2 \}$. 
Also, 
$\mu (C_{2, j}) = \mu (B_{2,k})$ for $j \in \{ 1,2, \ldots , p_3 \}$ 
and $k \in \{ 1,2, \ldots , q_3 \}$. 
Let $I_1 = B_{1,1}$ and $I_2 = B_{2,1}$. 
Stack the sets $B_{1,j}$ and $C_{1,k}$ such that whenever the sum 
of the values is negative, place a $B$ next, and otherwise place a $C$ set next. 
Stack the sets $B_{2,j}$ and $C_{2,k}$ such that whenever the sum 
of the values is negative, place a $B$ next, and otherwise place a $C$ set next. 
As long as $\delta_1 <  { \min{ \{ a, b \} } }/{2}$, then we have the precise number 
of level sets $B$ and $C$ to complete the two towers. 

For $x \in I_1$, 
\[
| \sum_{i=0}^{h_1 - 1} f( T^i x ) | = | p_2 b - q_2 c | = \delta_2 < \epsilon . 
\]
and 
$x \in I_2$, 
\[
| \sum_{i=0}^{h_2 - 1} f( T^i x ) | = | p_3 b - q_3 c | = \delta_1 - \delta_2 < \epsilon . 
\]
Equation (\ref{nbig4}) holds due to the greedy stacking algorithm used. 
The other conditions in the lemma hold by construction. $\Box$
\end{pf}

\begin{lem}
\label{dec1}
Suppose $f : X \to \real$ is a mean-zero finite step function. 
In particular, let $f = \sum_{i=1}^{m} a_i I_{I_i}$ where 
$X = \cup_{i=1}^{m} I_i$ is a disjoint union and $a_i$ are distinct 
real numbers for $1\leq i \leq m$ and $m \geq 2$. 
There exist disjoint measurable sets $J_1, J_2, \ldots , J_{m-1}$ such that 
$f$ takes on at most two values a.e. on $J_i$ and 
$\int_{J_i} f d\mu = 0$ for $1\leq i \leq m-1$. 
\end{lem}
\begin{pf}
We prove this by induction on $m$.  Clearly, this is true 
for $m=2$.  Suppose it is true for $m=n$. 
Let $m = n + 1$.  Choose $j$ such that for $1\leq i \leq n+1$, 
\[
\int_{I_j} | f | d\mu = | a_j | \mu ( I_j ) \leq \int_{I_i} | f | d\mu = | a_i | \mu ( I_i ) . 
\]
If $a_j \leq 0$, choose $k \neq j$ such that $a_k \geq 0$, otherwise choose 
$k$ such that $a_k \leq 0$.  
Choose $I^{\prime} \subset I_k$ such that 
\[
a_j \mu ( I_j ) + a_k \mu ( I^{\prime} ) = 0 . 
\]
Define $J_n = I_j \cup I^{\prime}$.  Thus, $f$ takes on at most $n$ steps 
on the subset $X \setminus J_n$. 
By induction, there exists $J_1, J_2, \ldots , J_{n-1}$ such that 
$f$ takes on at most two steps on $J_i$.  Therefore, our lemma is proved by induction. 
$\Box$
\end{pf}

\begin{lem}
\label{dec2}
Suppose $f: X \to \real$ is bounded and mean-zero.  Given $\epsilon_i > 0$ for $i\in \natural$, 
there exist a measure preserving map $T$, disjoint sets $I_i \subset X$ and natural numbers $h_i$ 
such that 
\begin{itemize}
\item $X = \bigcup_{i=1}^{\infty} \bigcup_{j=0}^{h_i-1} T^j I_i$ is a disjoint union, 
\item $ | \sum_{j=0}^{h_i-1} f ( T^j x ) | < \epsilon_i$ for $x\in I_i$, and 
\item $  | \sum_{j=0}^{k} f ( T^j x ) | < || f ||_{\infty} + \epsilon_i$ for $x\in I_i$ and $0 \leq k < h_i$. 
\end{itemize}
\end{lem}
\begin{pf}
If $f$ is a finite step function, then the lemma follows by applying Lemmas \ref{dec1} 
and \ref{two-step-lem} with a finite number of sets $J_i$.  
If $f$ is not a finite step function, then we apply Lemma \ref{ucl}, iteratively and potentially 
infinitely many times, to construct a sequence of TUB towers that satisfy this lemma. $\Box$
\end{pf}

\subsection{Proof of the Main Positive Result}
Now we are ready to proceed with the proof of Theorem \ref{cet}.
\medskip

\noindent \begin{pf}
Without loss of generality, we prove this theorem for the case $X = [0,1)$ and
$\mu$ equal to Lebesgue measure.  Also, we may assume $f \notin L^{\infty}$, 
since this case was handled previously, \cite{TAJR} and proposition \ref{bnd-prop}. 
If $f$ does not take on essentially infinitely many bounded values on a compact set, 
then first apply Lemma \ref{two-step-lem} to generate countable towers and transformation 
such that the sums are bounded, i.e., less than $\epsilon_i$ for $i^{th}$ tower. 
Let $k$ be the minimum positive integer such that $\mu( \{ x: 0 < f(x) \leq k \} ) > 0$,
and similarly let $\ell$ be the minimum positive integer such that $\mu( \{ x: 0 > f(x) \geq -\ell \} ) > 0$.
If no such $k$ and no such $\ell$ exist, then $f$ must equal zero almost everywhere,
and there is nothing left to prove.
Let $X_1 = \{ x: k-1 < f(x) \leq k \} $ and $Y_1 = \{ x: 1 - \ell > f(x) \geq -\ell \}$.
If $\int_{X_1} f d\mu + \int_{Y_1} f d\mu \leq 0$, define $Y_1^{\prime} \subseteq Y_1$
such that
\[
\int_{X_1} f d\mu + \int_{Y_1^{\prime}} f d\mu = 0 .
\]
In this case, let $X_1^{\prime} = X_1$.
Otherwise, choose $X_1^{\prime} \subset X_1$ such that
\[
\int_{X_1^{\prime}} f d\mu + \int_{Y_1} f d\mu = 0 .
\]
In this case, set $Y_1^{\prime} = Y_1$.  Also, define $k_1 = k$, $\ell_1 = \ell$ and
$X_0 = \{ x: f(x) = 0 \}$.
We may continue this procedure inductively to choose disjoint sets
$X_i^{\prime}, Y_i^{\prime}$ for $i=1, 2, \ldots$, and sequences
of positive integers $k_i, \ell_i$ such that
\begin{enumerate}
\item $k_i - 1 < f(x) \leq k_i$ for $x \in X_i^{\prime}$;
\item $1 - \ell_i  > f(x) \geq - \ell_i$ for $x \in Y_i^{\prime}$;
\item $k_{i+1} \geq k_i$, $\ell_{i+1} \geq \ell_i$;
\item $\lim_{i\to \infty} k_i + \ell_i = \infty$;
\item $\int_{X_i^{\prime}} f d\mu + \int_{Y_i^{\prime}} f d\mu = 0$;
\item $\mu \big( \bigcup_{i=1}^{\infty} ( X_i^{\prime} \cup Y_i^{\prime} ) \big) = \mu (X\setminus X_0)$.
\end{enumerate}
Let $\epsilon_i > 0$ such that $\sum_{i=1}^{\infty} \epsilon_i < \infty$. 
Let $\delta_i > 0$ such that 
\[
\delta_i ( \max{ \{ k_i, \ell_i \} } + \epsilon_i ) < \epsilon_i . 
\]
For each $i\in \natural$, apply Lemma \ref{dec2} to $f$ defined
on $X_i^{\prime} \cup Y_i^{\prime}$ to obtain a decomposition 
into potentially infinitely many towers satisfying the conditions 
of Lemma \ref{dec2}. 
Since the function $f_A$ which sums $f$ from the bottom to the top of each tower is bounded, 
then we can apply Proposition \ref{bnd-prop} to construct 
an ergodic measure preserving transformation $T_A$ on the base of the towers (called $A$) 
and bounded transfer function $g_A$ such that 
$f_A$ is a coboundary for $T_A$ with transfer function $g_A$.  
Then apply Lemma \ref{cob-ext-lem} to show that the ergodic measure preserving 
transformation $T$ on $X$ has coboundary $f$ with transfer function $g$ such that 
$| g(x) | < \max{ \{ k_i, \ell_i \} } + \epsilon_i$ 
for $x \in (X_i^{\prime} \cup Y_1^{\prime})$. 
For $i \in \natural$, let $Z_i = (X_i^{\prime} \cup Y_1^{\prime})$.
Since we assume that $f \notin L^{\infty}$, 
either $k_i \to \infty$ or $\ell_i \to \infty$.  Assume without loss of generality
that $k_i \to \infty$ and $k_i > \ell_i$ for infinitely many $i \in \natural$.
Either there exists $j$ such that $\ell_j > 1$, or we may choose $\ell_i \in (0, 1]$ 
such that $\ell_i$ is rational, nondecreasing, and
$Y_1^{\prime} = \{ x\in X: 0 > f(x) \geq -\ell_1 \}$ and
$Y_i^{\prime} = \{ x\in X: -\ell_{i-1} > f(x) \geq -\ell_i \}$
for $i \geq 2$. 
In either case, there exists $j$ such that 
$k_j \geq \ell_j$, and for $i > j$, $k_j + \epsilon_j \leq 2 ( k_j - 1 )$, 
and also for $x\in Y_{i}^{\prime}$, 
\[
-\ell_{i+1} \leq f(x) < -\ell_i < 0 .
\]
Thus,
\begin{eqnarray}
\int_{X} | g |^{p - 1} d\mu &=&
\sum_{i=1}^{\infty} \int_{Z_i} | g |^{p-1} d\mu \\
&=& \sum_{i=1}^{j} \int_{Z_i} | g |^{p-1} d\mu +
\sum_{i=j+1}^{\infty} \int_{Z_i} | g |^{p-1} d\mu
\end{eqnarray}
Since $|g|^{p-1}$ is bounded on $Z_i$ for $i \leq j$,
then $\sum_{i=1}^{j} \int_{Z_i} | g |^{p-1} d\mu < \infty$.
For each $i \in \natural$, let $m_i = \max{ \{ k_i, \ell_i \} }$. 
Also, if $m_i = k_i$, let $V_i = X_i^{\prime}$ and 
$W_i = Y_i^{\prime}$. Oherwise, let $V_i = Y_i^{\prime}$ and 
$W_i = X_i^{\prime}$. 
Thus, 
\begin{eqnarray}
\int_{V_i} | f | d\mu &\leq& 
\int_{V_i} (m_i) d\mu \\ 
&=& (m_i) \mu ( V_i ) \\ 
&=& \frac{(m_i)\mu(V_i)}{(\ell_{j})\mu(W_i)} ( \ell_j ) \mu ( W_i ) \\ 
&\leq& \frac{(m_i)\mu(V_i)}{(\ell_{j})\mu(W_i)} \int_{W_i} | f | d\mu 
\end{eqnarray} 
This implies
\[
\frac{ \mu( W_i ) }{ m_{i} } \leq
\frac{ \mu(V_i) }{ \ell_{j} } .
\]
Thus,
\begin{eqnarray*}
\sum_{i=j+1}^{\infty} \int_{Z_i} |g|^{p-1} d\mu
&=& \sum_{i=j+1}^{\infty} \int_{V_i} |g|^{p-1} d\mu +
\sum_{i=j+1}^{\infty} \int_{W_i} |g|^{p-1} d\mu \\
&\leq& \sum_{i=j+1}^{\infty} \int_{V_i} (m_i + \epsilon_i)^{p-1} d\mu +
\sum_{i=j+1}^{\infty} \int_{W_i} (m_i + \epsilon_i)^{p-1} d\mu \\
&\leq& \sum_{i=j+1}^{\infty} \int_{V_i} (m_i + \epsilon_i)^{p-1} d\mu +
\sum_{i=j+1}^{\infty} \int_{V_i} (m_i + \epsilon_i)^{p} \frac{1}{\ell_j}d\mu \\ 
&=& \sum_{i=j+1}^{\infty} \int_{V_i} (m_i - 1)^{p-1} \frac{(m_i + \epsilon_i)^{p-1}}{(m_i - 1)^{p-1}} d\mu +
\frac{1}{\ell_j} \sum_{i=j+1}^{\infty} \int_{V_i} (m_i - 1)^{p} \frac{(m_i + \epsilon_i)^{p}}{(m_i - 1)^{p}} d\mu \\ 
&\leq& 2^{p-1} \sum_{i=j+1}^{\infty} \int_{V_i} (m_i - 1)^{p-1} d\mu +
\frac{2^p}{\ell_j} \sum_{i=j+1}^{\infty} \int_{V_i} (m_i - 1)^{p} d\mu \\
&\leq& 2^{p-1} \sum_{i=j+1}^{\infty} \int_{V_i} | f |^{p-1} d\mu +
\frac{2^p}{\ell_j} \sum_{i=j+1}^{\infty} \int_{V_i} | f |^{p} d\mu \\
&\leq& 2^{p-1}  \| f \|_{p-1}^{p-1} + \frac{2^p}{\ell_j} \| f \|_p^{p} < \infty . 
\end{eqnarray*}
This completes the proof that $g \in L^{p-1}(X)$. $\Box$
\end{pf}

\section {\bf Non-existence of $L^p$-coboundaries}\label{cex}
In \cite{Ko2002}, Kornfeld shows that given $T \in \mathcal{E}$ 
which is a homeomorphism on a compact space $X$, 
there exists a continuous and bounded coboundary $f$ such that its 
associated transfer function is measurable, but not integrable. 
Also, it is pointed out that given $T$, $f$ may be constructed such that 
the transfer function $g$ is in $L^p$ for specified $p \geq 1$, 
but not contained in $L^q$ for $q > p$. 
However, if the function $f\in L^p_0$ is specified first, 
Kornfeld conjectured that there always exist 
an ergodic invertible measure preserving transformation $T$ and $g \in L^p$ 
such that $f = g - g \circ T$ a.e.\footnote{Kornfeld conveyed this conjecture to the second 
author verbally or through email.} 
In this section, we disprove this conjecture.  Furthermore, we prove 
a strong non-existence result showing that for generic $f \in L^p_0$, there are 
no $T \in \mathcal{E}$ and $g \in L^q$ for $q > p - 1$ such that 
$f = g - g \circ T$ a.e.  This is the statement of Theorem \ref{neg-thm}, 
and shows that generic $L^p_0$ functions lead to "wild" transfer functions 
(as termed in \cite{Ko2002}), 
universally for all $T \in \mathcal{E}$.  
Remark 2 in \cite{Kw} provides an argument for the existence of $L^{p}_0$ 
functions $f$ for $p \geq 2$ which are not coboundaries for any ergodic measure 
preserving transformation $T$ with transfer function $g \in L^{p}$.  
The argument in \cite{Kw} can be extended to show there are functions 
$f \in L^{p}_0$ which are not coboundaries for any ergodic measure 
preserving transformation $T$ with transfer function $g \in L^{q}$ 
for $q > p - 1$.   This is shown at the end of this section. 
Our results show this situation is generic for $f \in L^{p}_0$. 

A principal obstacle to solving the coboundary equation is imbalance between the positive 
and negative parts of a typical function $f \in L^p_0$. 
Suppose $a_i \in \real$ for $i\in \natural$ is an increasing sequence of real numbers 
such that $\lim_{i\to \infty} a_i = \infty$,
and for all reals $\alpha > 0$,
\begin{eqnarray}
\lim_{i\to \infty} \frac{a_i}{a_{i+1}^{\alpha}} = 0 . \label{cond1} 
\end{eqnarray}
Given $f\in L^p$ and $i \in \natural$, let 
\[
u_i(f) = \{ x \in X : f(x) < - a_i \} \ \mbox{and}\ v_i(f) = \{ x \in X : f(x) > a_i \} . 
\]
We are ready to define our generic class of $L^p_0$ functions. 
Given $n \in \natural$, define 
\[
\mathcal{G}^p_n = \{ f\in L^p_0 : \exists i > n\ \mid \ \mu(v_i(f)) > \frac{1}{a_i^p i^2}\ 
\mbox{and}\ \mu(u_{i-1} (f)) < \frac{1}{a_{i+1}^{p} i^2} \} . 
\]
Below we prove that $\mathcal{G}^p_n$ is both open and dense, and 
$f \in \cap_{n=1}^{\infty} \mathcal{G}^p_n$ satisfies the required property. 
The key property of the sequence $a_n$ is the fast growth rate. 
The following lemma will be used to guarantee that coboundaries 
$f\in \cap_{n=1}^{\infty} \mathcal{G}^p_n$ 
do not have transfer functions in $L^q$ for $q > p - 1$. 
\begin{lem}
\label{fast-lem}
For any $\alpha > 0$, 
\[
\lim_{n\to \infty} \frac{a_{n+1}^{\alpha}}{a_n n^2} = \infty . 
\]
\end{lem}
Now we prove that $\mathcal{G}^p_n$ is dense in $L_p$ for each $p\geq 1$ 
and $n\in \natural$. 
\begin{lem}
\label{lem-dense}
For each $n \in \natural$, the set $\mathcal{G}^p_n$ is dense in $L^p_0$. 
\end{lem}
\begin{pf}
Let $f \in L^p_0$, $\epsilon > 0$ and $n\in \natural$. 
Since bounded measurable functions are dense in $L^p_0$, 
we can choose a bounded mean zero $f_0 \in L^p$ 
such that 
\[
|| f - f_0 || < \frac{\epsilon}{3} . 
\]
Choose $i_0 \geq n$ such that $a_{i_0} > || f_0 ||_{\infty}$. 
Choose $i_1 \geq i_0$ such that 
\[
\frac{2^{p+2}}{i_1^2} < \frac{\epsilon}{3} . 
\]
Choose a subset $Y \subset X$ such that 
\[
\mu(Y) = \frac{2}{a_{i_1}^{p} i_1^2 } + \frac{4}{a_{i_1-1}a_{i_1}^{p-1} i_1^2 } 
\]
and $\int_Y f d\mu = 0$. 
Let $V\subset Y$ be such that $\mu(V) = 2 \big( a_{i_1}^p i_1^2 \big)^{-1}$ 
and define $U = Y\setminus V$. 
Define $f_1$ as a modification of $f_0$ in the following manner:
\begin{eqnarray*} 
f_1(x)= 
\left\{\begin{array}{ll}
2 a_{i_1}, & \mbox{if $x\in V$}, \\ 
- a_{i_1 - 1} & \mbox{if $x\in U$}, \\ 
f_0(x), & \mbox{if $x\in X\setminus Y$} . 
\end{array}
\right.
\end{eqnarray*}
Thus,
\begin{align*}
|| f - f_1 ||_{p} &\leq || f - f_0 ||_{p} + || f_0 - f_1 ||_{p} \\ 
& < \frac{\epsilon}{3} + 2^p a_{i_1}^{p} \mu(V) + a_{i_1-1}^{p} \mu(U) \\ 
& < \epsilon . 
\end{align*}
Also, $f_1 \in \mathcal{G}^p_n$ which completes the proof. $\Box$ 
\end{pf}

\begin{lem}
\label{lem-open}
For each $n \in \natural$, the set $\mathcal{G}^p_n$ is open in $L^p_0$. 
\end{lem}
\begin{pf}
Suppose $f \in \mathcal{G}^p_n$.  Then there exists $i \geq n$ such that 
\begin{align*}
\mu_i &= \mu(v_i(f)) > \frac{1}{a_i^p i^2} ,\\ 
\nu_i &= \mu(u_{i-1}(f)) < \frac{1}{a_{i+1}^{p} i^2} . 
\end{align*}
Thus, there exists $a^{\prime} > a_i$ and $a^{\prime \prime} > a_{i-1}$, 
and $\mu^{\prime}, \nu^{\prime}$ such that 
\begin{align*}
\mu( \{ x: f(x) > a^{\prime} \} ) &> \mu^{\prime} > \frac{1}{a_i^p i^2 } , \\ 
\mu( \{ x: f(x) < - a^{\prime \prime} \} ) &< \nu^{\prime} < 
\frac{1}{a_{i+1}^{p} i^2 } . 
\end{align*}
Define $\epsilon > 0$ as 
\[
\epsilon = \min{ \{ ( \mu^{\prime} - \frac{1}{a_i^p i^2} ) ( a^{\prime} - a_i )^{p} , 
( \frac{1}{a_{i+1}^{p} i^2 } - \nu^{\prime} ) (a^{\prime \prime} - a_{i-1})^{p} \} } . 
\]
It is not difficult to see that the $\epsilon$-ball centered at $f\in L^p_0$ is contained 
in $\mathcal{G}^p_n$. 
$\Box$
\end{pf}

\noindent 
Let 
\[
\mathcal{G}_p = \bigcap_{n=1}^{\infty} \mathcal{G}^p_n . 
\]
We have the following core result of this paper. 

\begin{prop}
\label{core-prop}
Suppose $f \in \mathcal{G}_p$, $T\in \mathcal{E}$ and $g$ is a measurable function.
If the coboundary equation $f = g - g\circ T$ is satisfied a.e., 
then $g \notin L^q$ for $q > p - 1$. 
\end{prop}

\noindent 


\noindent 
{\bf Proof of Proposition \ref{core-prop} }: 
Let $\sgn$ be the standard sign function defined as
$\sgn(i)=-1$, if $i<0$, $\sgn(i)=0$, if $i=0$ and
$\sgn(i)=1$, if $i>0$.
For $i \in \integers$, let $[ i ] = \{ j \in \integers : i \leq j < 0 \}$ if $i < 0$,
and $[ i ] = \{ j \in \integers : 0 \leq j < i \}$ if $i \geq 0$.
Note, for $i \in \integers$, the coboundary equation expands to the following:
\[
g(T^i x) = g(x) - \sgn (i) \sum_{j \in [i]} f ( T^j x) .
\]
Define our specialized sign function $\rho: X \to \{ -1, 1 \}$ based on the following:
\begin{enumerate}
\item if $g(x) \leq {a_n}/{2}$, let $\rho(x)=1$, \label{prop-case-1}
\item otherwise if $g(x) > {a_n}/{2}$, then 
let $\rho(x)=-1$. \label{prop-case-2}
\end{enumerate}
For $n \in \natural$, let $c_n = {a_{n}} / {a_{n-1}}$. 
Assume $f \in L^p(X)$ and $q > p - 1$.
Choose integer $k > 1$ such that $kq > p$. 
Let $A_n = u_{n-1}(f)$ and $B_n = v_{n}(f)$. 
For $x \in B_n$, let 
\[
\ell_x = \min{ \{ \ell : c_n 
\leq \ell < \lceil ( c_n )^{k} \rceil ,
| g( T^{ \rho(x) \ell } x ) | < \frac{ a_{n} }{4} \big( c_n \big)^{h-1} ,
\lceil ( c_n )^{h} \rceil \leq \ell < \lceil ( c_n )^{h+1} \rceil
\} },
\]
otherwise, let $\ell_x = (c_n)^k$. 
Given $x \in X$, define the set $L(x) = [ \rho(x) \ell_x ]$. 
Choose $N\in \natural$ such that for $n \geq N$, 
\[
\frac{a_{n+1}^{p}}{a_{n}^{k+p} }  > 4 . 
\]
Thus, for $n \geq N$,
\[
\mu \Big( \bigcup_{j=-c_n^{k}}^{c_n^{k}} T^{j} (A_n ) \Big) <  
\frac{2 a_n^k}{a_{n+1}^p n^2} < \frac{1}{2} \frac{1}{a_n^p n^2} < \frac{1}{2} \mu (B_n) .
\]
Let
\[
B_n^{\prime} = B_n \setminus
\Big( \bigcup_{j=-c_n^{k}}^{c_n^{k}} T^{j} (A_n ) \Big) .
\]
Hence, $\mu (B_n^{\prime}) > \frac{1}{2} \mu (B_n)$ for $n \geq N$.
We break the proof down into 4 separate cases and handle
each separately.
\medskip

\begin{enumerate}
\item $B_{n,1} = \{ x \in B_n^{\prime} : g(x) < {a_n}/{2}, \ell_x < \lceil ( c_n )^{k} \rceil \}$,
\item $B_{n,2} = \{ x \in B_n^{\prime} : g(x) < {a_n}/{2}, \ell_x \geq \lceil ( c_n )^{k} \rceil \}$,
\item $B_{n,3} = \{ x \in B_n^{\prime} : g(x) \geq {a_n}/{2}, \ell_x < \lceil ( c_n )^{k} \rceil \}$,
\item $B_{n,4} = \{ x \in B_n^{\prime} : g(x) \geq {a_n}/{2}, \ell_x \geq \lceil ( c_n )^{k} \rceil \}$.
\end{enumerate}
\medskip

At least one of the $B_{n,m}$ satisfies $\mu(B_{n,m}) \geq ({1}/{8}) \mu(B_n)$
for $m=1,2,3,4$.  We handle the case $\mu(B_{n,1}) \geq ({1}/{8}) \mu(B_n)$ first.
We create tiles in the following way. 
For $x \in B_{n,1}$, let
\[
B_{n,1}(x) = \{ T^i x : i \in L(x) \} .
\]
There exists $J = J_{n,1}$ such that for $x \neq y$, $x, y \in J$,
\[
B_{n,1}(x) \cap B_{n,1}(y) = \emptyset
\]
and
\[
\mu ( B_{n,1} \cap \bigcup_{x \in J} B_{n,1}(x) ) > \frac{1}{2} \mu (B_{n,1}) .
\]

The $L^q$-norm of the transfer function $g$ will blow up on the set $J$. 
Before completing the general proof, it is helpful to see
how the argument goes in a special case.
Suppose $\ell_x = \frac{a_n}{a_{n-1}}$ for $x \in J$.
This implies for $x\in J$, $T^{i}(x) \notin B_n$ on the order
of $\frac{a_n}{a_{n-1}}$ times.  Also, for this special case,
$T^i (x)$ cannot fall in $B_n$ for $0 < i < \ell_x$.
Note that $T^i(x)$, $0\leq i < \ell_x$, does not fall in $A_n$
by the previous choice of $J$.
But, for $x\in J$, the transfer function at $T^i(x)$
will be on the order of the sum, so that $g(T^i(x))$ will
be on the order of $a_n$ (or $\frac{a_n}{4}$).
This implies
\begin{align}
\int_{X} | g(x) |^q d\mu &\approx
\Big( \frac{a_n}{4} \Big)^{q} \Big( \frac{a_n}{a_{n-1}} \Big) \mu(B_n) \\
&= \frac{1}{4^q} \frac{a_n^{q+1}}{a_{n-1} a_n^p n^2} \\
&= \frac{1}{4^q} \frac{a_n^{q+1-p}}{a_{n-1} n^2}
\end{align}
But the last term tends to infinity as $n\to \infty$ by
the definition of $a_n$ and Lemma \ref{fast-lem}. 

\noindent
{\bf General proof for case 1}:
First we prove the following lemma.
\begin{lem}
Suppose $\lceil (c_n)^h \rceil \leq \ell_x < \lceil (c_n)^{h+1} \rceil$
for $1 \leq h < k$.  If
\[
\ell_0 = \# \{ i \in L(x) : T^i(x) \in B_n \} ,
\]
then
\[
\ell_x > \frac{1}{2} c_n \ell_0 .
\]
\end{lem}

\noindent
{\bf Proof of lemma}:
Suppose the lemma is not true.
Then 
\begin{align}
| g(x) - \sgn (x) \sum_{i\in L(x) } f(T^i x) |
&\geq \ell_0 a_n - ( \ell_x - \ell_0 ) a_{n-1} = 
\ell_0 a_n - \ell_x a_{n-1} + \ell_0 a_{n-1} \\
&\geq a_{n-1} \ell_x + \frac{2a_{n-1}^2}{a_n} \ell_x \\
&\geq \frac{a_n^h}{a_{n-1}^{h-1}} + 2\frac{a_n^{h-1}}{a_{n-1}^{h-2}}
> a_n (c_n)^{h-1} .
\end{align}
This contradicts the definition of $\ell_x$. $\Box$

\noindent
{\bf Resume proof of proposition}:
\medskip

Thus, we have the following 
\begin{align}
\int_{X} | g(x) |^q d\mu &\geq
\int_{J} \sum_{i \in L(x)} | g(T^i x) |^q d\mu \\
&= \int_{J} \sum_{i \in L(x)} | g(x) - \sum_{j \in [i] } f(T^j x) |^q d\mu \\
&\geq \int_{J} \sum_{i \in L(x)} | \frac{a_{n}}{4} |^q d\mu \\
&= | \frac{a_n}{4} |^q \ell_x \mu (J)  d\mu \\
&> | \frac{a_n}{4} |^q \int_{J} \frac{1}{2} \big( \frac{a_n}{a_{n-1}} \big)
\sum_{i \in L(x) } I_{B_n} (T^i x) d\mu \\
&>  | \frac{a_n}{4} |^q \frac{1}{2} \big( \frac{a_n}{a_{n-1}} \big)
\big( \frac{1}{2} \mu (B_{n,1}) \big) \\
&>  \Big( \frac{1}{32} \Big) | \frac{a_n}{4} |^q \big( \frac{a_n}{a_{n-1}} \big) \mu (B_{n}) \\
&= \frac{a_n^{q+1}}{32 (4^q) a_{n-1} a_n^p n^2} \\
&= \frac{a_n^{q+1-p}}{32 (4^q) a_{n-1} n^2} . 
\end{align}
The proof for this case is complete, since, by condition (\ref{cond1}),
\[
\lim_{n\to \infty} \frac{a_n^{q+1-p}}{32 (4^q) a_{n-1} n^2} = \infty .
\]

\noindent
{\bf Proof for case 2}:
\begin{align}
\int_{B_n} | g(x) |^q d\mu &\geq \int_{J} \sum_{i \in L(x)} | g(T^i x) |^q I_{B_n} (T^i x) d\mu \\
&= \int_{J} \sum_{i \in L(x)}  | g(x) - \sgn(x) \sum_{ j \in [-i] } f(T^j x) |^q I_{B_n} (T^i x) d\mu \\
&\geq \int_{J} \sum_{i = c_n^{k-1}}^{c_n^k}  | g(x) - \sgn(x) \sum_{ j \in [-i] } f(T^j x) |^q I_{B_n} (T^{-i} x) d\mu \\
&\geq \sum_{i = c_n^{k-1}}^{c_n^k} \int_{J} \Big( \big( \frac{a_n}{4} \big)
\big( \frac{a_n}{a_{n-1}} \big)^{k-1} \Big)^q I_{B_n} (T^i x) d\mu\\
&\geq \sum_{i = c_n^{k-1}}^{c_n^k} \int_{J} \Big( \big( \frac{a_n}{4} \big)
\big( \frac{a_n}{a_{n-1}} \big)^{k-1} \Big)^q I_{B_n} (T^i x) d\mu \\
&\geq \Big( \frac{c_n^k - c_n^{k-1}}{c_n^k} \Big) \big( \frac{1}{16} \big) \mu(B_n)
\Big( \big( \frac{a_n}{4} \big) \big( \frac{a_n}{a_{n-1}} \big)^{k-1} \Big)^q \\
&> \big( \frac{1}{32} \big) \Big( \frac{a_n^{kq}}{ a_n^p n^2 4^q a_{n-1}^{q(k-1)} } \Big) \\
&= \big( \frac{1}{32} \big) \Big( \frac{a_n^{kq-p}}{ n^2 4^q a_{n-1}^{q(k-1)} } \Big) 
\end{align}
Since 
\[
\lim_{n\to \infty} \big( \frac{1}{32} \big) \Big( \frac{a_n^{kq-p}}{ n^2 4^q a_{n-1}^{q(k-1)} } \Big) 
= \infty ,
\]
then our result follows for case 2. 

Case 3 would be handled in a similar manner as case 1, 
except we would base our estimate of $g(x)$ on the inverse of $T$. 
Thus, we have the following 
\begin{align}
\int_{X} | g(x) |^q d\mu &\geq
\int_{J} \sum_{i \in L(x)} | g(T^i x) |^q d\mu \\
&= \int_{J} \sum_{i \in L(x)} | g(x) + \sum_{j \in [i] } f(T^j x) |^q d\mu \\
&\geq \int_{J} \sum_{i \in L(x)} | \frac{a_{n}}{4} |^q d\mu 
\end{align}
The next steps continue in a similar manner as case 1. 
Also, case 4 follows in a similar manner as case 2, 
except by using $T^{-1}$ instead of $T$. 
$\Box$


\noindent 
{\bf Proof of Theorem \ref{neg-thm}}: 
Define 
\[
\mathcal{G}_p = \bigcap_{n=1}^{\infty} \mathcal{G}^p_n . 
\]
By Lemmas \ref{lem-dense} and \ref{lem-open}, 
the set $\mathcal{G}_p$ is a dense $G_{\delta}$ subset of $L^p_0$. 
Also, by Proposition \ref{core-prop}, $f \in \mathcal{G}_p$ satisfies 
the conditions of Theorem \ref{neg-thm}. 
$\Box$ 

\subsection {\bf Not a moment}\label{cex2}
Let $\phi: \real \to \real$ be a measurable function such that
\[
\lim_{x \to \infty} \phi (x) = \infty .
\]
For $i \in \natural$, let $A, B_i$ be disjoint sets in $X$, and $b_i > 0$.  Define $f$ as
\[
f = I_{A} - \sum_{i=1}^{\infty} b_i I_{B_i} .
\]
We will give conditions on the fast growth rate of $b_i$ as well as conditions on the sets
$A, B_i$ to guarantee that $f$ is contained in $L^1$, but such that
$\phi \circ \left\vert g \right\vert$ is not in $L^1$ for any transfer function $g$ 
of an ergodic invertible measure preserving transformation $T$.
Let $A \subset X$ have measure $\mu (A) = {1}/{2}$.
Choose $b_i > 0$ for $i\in \natural$ such that
$\lim_{i\to \infty} b_i = \infty$,
and such that for all reals $\alpha > 0$, 
\begin{eqnarray}
\lim_{i\to \infty} \frac{b_i}{b_{i+1}^{\alpha}} = 0 \label{cond6}
\end{eqnarray}
and also for $y \geq {b_i} / {4}$, 
\begin{eqnarray}
\frac{\phi ( y )}{2^i} \geq i . \label{cond7}
\end{eqnarray}
Choose disjoint sets $B_i \subset A^c$ such that
\begin{eqnarray}
\mu (B_i) = \frac{1}{b_i 2^{i+1}}.\label{cond9}
\end{eqnarray}
Observe that $f\in L^1$ is mean zero.
\begin{prop}
\label{cex-prop2}
Let $\phi: \real \to \real$ be a measurable function satisfying
$\lim_{x\to \infty} \phi (x) = \infty$.
Suppose the mean zero function
$f = I_{A} - \sum_{i=1}^{\infty} b_i I_{B_i}$ satisfies the conditions above,
including (\ref{cond6}, \ref{cond7}, \ref{cond9}).
If $T$ is an ergodic invertible measure preserving transformation $T: X \to X$
and $g$ is a transfer function satisfying
$f(x) = g(Tx) - g(x)$ for almost every $x\in X$, then
\[
\int_X \phi ( \left\vert g \right\vert ) d\mu = \infty .
\]
\end{prop}

\noindent
\begin{pf}
Let $\sgn$ be the standard sign function defined as
$\sgn(i)=-1$, if $i<0$, $\sgn(i)=0$, if $i=0$ and
$\sgn(i)=1$, if $i>0$.
For $i \in \integers$, let $[ i ] = \{ j \in \integers : i \leq j < 0 \}$ if $i < 0$,
and $[ i ] = \{ j \in \integers : 0 \leq j < i \}$ if $i \geq 0$.
Note, for $i \in \integers$, the coboundary equation expands to the following:
\[
g(T^i x) = g(x) + \sgn (i) \sum_{j \in [i]} f ( T^j x) .
\]
Define our specialized sign function $\rho: X \to \{ -1, 1 \}$ based on the following:
\begin{enumerate}
\item if $g(x) \leq {b_n}/{2}$, let $\rho(x)=1$, \label{prop-case-11} 
\item otherwise if $g(x) > {b_n}/{2}$, then
let $\rho(x)=-1$. \label{prop-case-22}
\end{enumerate}
Assume $f \in L^1$. 
For $x \in B_n$, let
\[
\ell_x = \min{ \{ \ell : \ell > 0 ,
| g( T^{ \rho(x) \ell } x ) | < \frac{ b_n }{4} \big( b_n \big)^{h-1} ,
\lceil ( b_n )^{h} \rceil \leq \ell < \lceil ( b_n )^{h+1} \rceil
\} },
\]
Note, $\ell_x < \infty$ for almost every $x \in X$, otherwise
our result follows directly.
Thus, exclude points $x \in X$ where $\ell_x = \infty$.
Choose $k_n \in \natural$ such that
\[
\mu ( \{ x \in B_n : \ell_x < b_n^{k_n+1} \} ) > \frac{1}{2} \mu (B_n) .
\]
Given $x \in X$, define the set $K_n(x) = [ \frac{1}{2} \rho(x) \ell_x ]$.
We do not need to consider all of the cases as in Proposition \ref{core-prop}, 
due to the special nature of the counterexamples $f$ in this result. 
We create tiles in the following way.
For $x \in B_n$, let
\[
B_{n}(x) = \{ T^i x : i \in K_n(x) \} .
\]
There exists $J_n$ such that for $x \neq y$, $x, y \in J_n$,
\[
B_{n}(x) \cap B_{n}(y) = \emptyset
\]
and
\[
\mu ( B_{n} \cap \bigcup_{x \in J_n} B_{n}(x) ) > \frac{1}{4} \mu (B_{n}) .
\]

\noindent
First we prove the following lemma.
\begin{lem}
Suppose $\lceil (b_n)^h \rceil \leq \ell_x < \lceil (b_n)^{h+1} \rceil$
for $1 \leq h < k_n+1$.  If
\[
\ell_0 = \# \{ i \in [\rho(x) \ell_x ] : T^i(x) \in B_n \} , 
\]
then
\[
\ell_x > \frac{1}{2} b_n \ell_0 .
\]
\end{lem}

\noindent
{\bf Proof of lemma}:
Suppose the lemma is not true.
Then
\begin{align}
| g(x) + \sgn (x) \sum_{i\in [\rho(x) \ell_x] } f(T^i x) |
&\geq \ell_0 b_n - ( \ell_x - \ell_0 ) =
\ell_0 b_n - \ell_x + \ell_0 \\
&\geq \ell_x + \frac{2}{b_n} \ell_x \\
&\geq b_n^h + 2 b_n^{h-1}
> b_n^{h} .
\end{align}
This contradicts the definition of $\ell_x$. $\Box$

\noindent
{\bf Resume proof of proposition}: 
Thus, we have the following
\begin{align}
\int_{X} \phi \Big( | g(x) | \Big) d\mu &\geq
\int_{J_n} \sum_{i \in K_n(x)} \phi \Big( | g(T^i x) | \Big) d\mu \\
&= \int_{J_n} \sum_{i \in K_n(x)} \phi \Big( | g(x) + \sum_{j \in [i] } f(T^j x) | \Big) d\mu \\
&\geq \int_{J_n} n 2^n \ell_x  d\mu \\
&> n 2^n \int_{J_n} \frac{1}{2} \big( b_n \big)
\sum_{i \in K_n(x) } I_{B_n} (T^i x) d\mu \\ 
&>  n 2^n \frac{1}{2} \big( b_n \big)
\big( \frac{1}{4} \mu (B_{n}) \big) \\
&>  \Big( \frac{1}{8} \Big) n 2^n \big( b_n \big) \mu (B_{n}) \\
&= \frac{ n 2^n }{8 (2^{n+1})} \to \infty , \mbox{ as } n\to \infty. \ \Box 
\end{align}
\end{pf}


\subsection{Kwapien argument for the non-existence of $L^{p}$ coboundaries}
The following proposition establishes the existence of $L^p$ functions $f$ 
with no transfer function in $L^q$ for $q > p - 1$.  
The argument is due to Kwapien \cite{Kw}. 
\begin{prop}[Remark 2 in \cite{Kw}]
Given $p \in \real$ such that $p \geq 2$, there exists $f \in L^{p}$ such that 
for any solution pair $T$ and $g$ to the equation $f = g - g\circ T$ where 
$T$ is an ergodic invertible measure preserving transformation, then 
$g \notin L^q$ for $q > p - 1$. 
\end{prop}

\begin{pf}
Let $f \in L^p$ be such that $\int f d\mu = 0$, $f(x) \geq -1$ for a.e. $x$, 
and for $r > p - 1$, 
\[
\limsup_{n \to \infty} \Big( n \int_{f > n} | f - n |^{r} \Big) d\mu = \infty . 
\]
We refer to this as the Kwapien condition. 
To obtain examples $f$ that satisfy the Kwapien condition, 
suppose $p \geq 2$, $r > p - 1$, and let $N_k \in \natural$ be such that 
\begin{eqnarray}
\label{eqn-kw2}
\sum_{k=1}^{\infty}  N_k^{ - {(1+r-p)}/{2(r+1)}  } < \frac{1}{2^{p+1}} . 
\end{eqnarray}
Let 
\[
\delta = \frac{ 1 + r - p }{ 2(r+1) } . 
\]
By (\ref{eqn-kw2}), 
\[
\lim_{k\to \infty} N_k^{ {(1+r-p)}/{2} } = \lim_{k\to \infty} N_k^{ \delta (r+1) } = \infty . 
\]
Let $E_k$ be disjoint sets for $k \in \natural$ such that 
$\mu (E_k) = {1}/{N_k^p}$.  Thus,  
\[
\sum_{k=1}^{\infty} N_k^{ - {p(1+r-p)}/{2(r+1)}  } < {1}/{2^{p+1}} . 
\]
Define $f^{+}$ such that 
\[
f^{+} = \sum_{k=1}^{\infty} 2 N_k^{ 1 - \delta } I_{E_k} . 
\]
We've shown that $\int f^{+} d\mu < {1} / {2}$. 
Let $E_0$ be a subset disjoint from $\bigcup_{k=1}^{\infty} E_k$ such that $\mu (E_0) = \int f^{+} d\mu$. 
Define $f = f^{+} - I_{E_0}$.  Thus, $\int f d\mu = 0$ and $|| f ||_p < \infty$. 

Let $L_k = n_k\int_{f > n_k} (f - n_k)^r d\mu$ where $n_k = N_{k}^{1-\delta}$.   
This $n_k$ is not a whole number probably, but I am going to ignore that.
Then $L_k \ge N_{k}^{1-\delta} \int_{E_{k}} ( N_{k_o}^{1-\delta})^r d\mu$.  We get 
\[
L_k \ge N_{k}^{(r+1)(1-\delta)}/N_{k_o}^p = N_{k}^{ { (1+r+p) } /{2} }/N_{k}^p = N_{k}^{(r+1 - p)/2} . 
\]
Since $\lim_{k\to \infty} L_k = \infty$ and $f$ satisfies the Kwapien condition. 

Now we prove that $f$ is not a coboundary with a transfer function in $L^r$ for 
any $r > p - 1$. 
Since $f \geq -1$ a.e., then for a.e. $x$, 
\[
\left| \sum_{i=0}^{n} f ( T^i x ) \right| \geq \sum_{i=0}^{n} \big( f ( T^i x ) - n \big)  \mathbb{I} \{ f (T^i (x)) > n \} . 
\]
Each term in the sum on the right side of the inequality is non-negative and 
therefore, 
\begin{eqnarray}
\left|\left| \sum_{i=0}^{n} f ( T^i x ) \right|\right|_r^r & \geq & 
\sum_{i=0}^{n} \int_{f\circ T^i > n} \big( f ( T^i x ) - n \big)^r  d\mu \\ 
& = & (n + 1) \int_{f > n} \big( f - n \big)^r d\mu . 
\end{eqnarray}
Therefore, 
$\sum_{i=0}^{n} f ( T^i x ) = g(x) - g ( T^{n+1} x )$ is unbounded in $L^r$. $\Box$ 
\end{pf}

\section{Category of Transformation Solutions}
\label{cat-sec}
In this section, we prove for any non-trivial measurable function $f$, 
the set of ergodic measure preserving transformations $T$ 
such that the coboundary equation $f = g - g\circ T$ has a measurable 
solution $g$ is a first category set (meager). 

\begin{prop}
Let $f$ be a measurable function such that $\mu (\{ x : f(x) \neq 0 \}) > 0$. 
Let $\mathcal{T}$ be the set of ergodic invertible measure preserving 
transformations $T$ such that $f = g - g\circ T$ has a measurable solution 
$g$.  The set $\mathcal{T}$ is a set of first category (meager). 
\end{prop}
\begin{pf}
Let $\eta \in \real$ such that $0 < \eta < {1} / {10}$. 
For each $n \in \natural$, define 
\[
D_n = \{ T \in \mathcal{E} : \exists k > n \ \mbox{such that}\ 
\mu ( \{ x: | \sum_{i=0}^{k-1} f(T^i x) | > n \} ) > \eta \} . 
\]
For each $n\in \natural$, the set $D_n$ is both open and dense. 
Establishing open-ness is straightforward. Let $T \in D_n$. 
There exists $\delta_0 > 0$ such that 
\[
\mu \big( \{ x: | \sum_{i=0}^{k-1} f(T^i x) | > n + \delta_0 \} \big)  >  \eta . 
\]
Thus, if $\eta_0 = \mu \big( \{ x: | \sum_{i=0}^{k-1} f(T^i x) | > n + \delta_0 \} \big)$, then 
\[
\big\{ S \in \mathcal{E} : \int_{X} | f\circ S^i - f\circ T^i | d\mu < \frac{\delta_0 ( \eta_0 - \eta )}{k} \big\} 
\]
is an open neighborhood containing $T$ and contained in $D_n$. 

To establish that $D_n$ is dense, it can be accomplished by an application of the ergodic theorem. 
Let $S\in \mathcal{E}$ and $\epsilon \in \real$ be such that ${1} / {20} > \epsilon > 0$.  
If $S \in D_n$, then we set $T = S$. 
Otherwise, assume $S \notin D_n$. 
Choose $\alpha > 0$ such that the set 
$A = \{ x \in X : f(x) > \alpha \}$ has positive measure. 
Similarly, choose $\beta > 0$ such that the set 
$B = \{ x \in X : f(x) < -\beta \}$ has positive measure. 
Let $\gamma \in \natural$ be such that 
\[
\gamma \geq \max{ \{ \frac{2n}{\alpha} , \frac{2n}{\beta} \} } . 
\]
Choose $\ell_0 > n$ such that for $\ell \geq \ell_0$, 
\[
\mu \big( \{ x \in X: \sum_{i=0}^{\ell-1} I_{A} ( S^{i} x ) > \gamma \} \big) > 1 - \epsilon 
\]
and 
\[
\mu \big( \{ x \in X: \sum_{i=0}^{\ell-1} I_{B} ( S^{i} x ) > \gamma \} \big) > 1 - \epsilon . 
\]
Choose $h > \ell_0$ such that 
\[
\frac{\ell_0}{h} < \frac{\epsilon}{4} . 
\]
There is a Rohklin tower of height $4h$ with base $I$ such that 
\[
\mu \Big( \bigcup_{i=0}^{4h-1} S^i I \Big) > 1 - \frac{\epsilon}{4h} . 
\] 
There exist disjoint sets $I_1, I_2 \subset I$ such that for each $x \in I_1$ and $y \in I_2$, 
there exists $j(x), j(y)$ such that $h \leq j(x) < 2h$, $h \leq j(y) < 2h$, and 
\[
\sum_{i=0}^{\ell_0 - 1} I_A (S^{i + j(x)} x) > \gamma 
\]
and 
\[
\sum_{i=0}^{\ell_0 - 1} I_B (S^{i + j(y)} y) > \gamma . 
\]
By the choice of $\epsilon < {1} / {20}$, then $I_1, I_2$ may be chosen 
such that 
\[
\mu (I_1) = \mu (I_2) > \frac{1}{4} \mu (I) . 
\]
For each $x \in I_1$, let $i_1(x), i_2(x), \ldots i_{\gamma}(x)$, be 
increasing such that 
\[
S^{ i_j(x) } (x) \in A 
\]
and similarly, 
for each $y \in I_2$, let $i_1(y), i_2(y), \ldots i_{\gamma}(y)$, be 
such that 
\[
S^{ i_j(y) } (y) \in B 
\]
and $h \leq i_j(x) < 2h - 1$, $i_{\gamma}(x) < i_{1}(x) + \ell_0$, 
and $h \leq i_j(y) < 2h - 1$, $i_{\gamma}(y) < i_{1}(y) + \ell_0$. 
Let $\phi: I_1 \to I_2$ be an invertible measure preserving map. 
The transformation $T$ will be defined in the following manner:
for $x \in I_1$, let $y = \phi(x) \in I_2$, 
\[
T^{i_j(x)} (x) = S^{i_j(y)} (y) 
\]
and 
\[
T^{i_j(y)} (y) = S^{i_j(x)} (x) . 
\]
Otherwise, define $T$ to be identical to $S$ everywhere else on $X$. 
Consider 
\[
\sum_{i=0}^{3h-1} f ( T^i x ) 
\]
for $x \in \bigcup_{i=0}^{h - 1} T^i (I_1 \cup I_2)$. 
Note for such $x$, 
\[
\Big| \sum_{i=0}^{3h-1} f ( T^i x ) - \sum_{i=0}^{3h-1} f ( S^i x ) \Big| > 2n . 
\]
Since 
\[ 
\mu \Big( \bigcup_{i=0}^{h - 1} T^i (I_1 \cup I_2) \Big) > \frac{1}{5} , 
\]
then 
\[
\mu \Big( \big\{ x \in X : \big| \sum_{i=0}^{3h-1} f ( T^i x ) \big| > n \big\} \Big) > \eta . 
\]
This implies $T \in D_n$ and $|| T - S || < \epsilon$.  
Thus, $\mathcal{T} = \cap_{n=1}^{\infty} D_n$ is a dense $G_{\delta}$ set. 
If $T \in \mathcal{T}$, then Schmidt's condition (\ref{schmidt}) 
for a measurable transfer function does not hold, and our result follows. $\Box$ 
\end{pf}

Since the set of $T \in \mathcal{E}$ that yield a measurable solution $g$ is meager, 
this raises the question of whether the set of solutions $T \in \mathcal{E}$, $g \in L^0$, 
is nonempty. 
In general, there are solutions, including cases where $f\notin L^1$. 
The following theorem extends the result of Anosov \cite{Anosov73} to show when 
there exists a solution pair $T \in \mathcal{E}$, $g \in L^0$. 
\begin{prop}
Suppose $f: X \to \real$ is a measurable function. 
The coboundary equation  $f = g - g \circ T$ has solutions $T \in \mathcal{E}$, 
$g \in L^0$, if and only if, 
\[
\int_{f > 0} f d\mu = \int_{f < 0} \big( - f \big) d \mu \ \ (\infty\mbox{ or finite}). 
\]
\end{prop}
\begin{pf}
The case where both $\int_{f > 0} f d\mu$ and $\int_{f < 0} \big( - f \big) d \mu$ 
are finite and unequal is already covered by Anosov's result \cite{Anosov73}. 
If the integrals are finite and equal, it follows from Theorem \ref{pos-thm}. 

Next, we prove the case where one integral is finite and the other is infinite. 
Without loss of generality, assume 
 $\int_{f > 0} f d\mu = \infty$ and $\int_{f < 0} \big( - f \big) d \mu < \infty$. 
Choose a measurable subset $A \subset \{ f > 0 \}$ such that 
\[
\int_{A} f d \mu + \int_{f < 0} f d \mu = 1 . 
\]
Let 
\begin{eqnarray*} 
f_0(x)= 
\left\{\begin{array}{ll}
f(x) & \mbox{if $x \in A \cup \{ f < 0 \} $}, \\ 
0 & \mbox{if $x \in \{ f > 0 \} \setminus A$}.
\end{array}
\right.
\end{eqnarray*}
Thus, $f_0 \in L^1$ and $\int_X f_0 d \mu = 1$. 
Given $T \in \mathcal{E}$, by the mean ergodic theorem, 
\[
\lim_{n \to \infty} \int_X \Big| 
\frac{1}{n} \sum_{i=0}^{n-1} f_0 (T^i x) - \int_X f_0 d \mu \Big| d \mu = 0 . 
\]
Let $\delta > 0$. 
Then 
$\mu \{ x \in X : \big| 
\frac{1}{n} \sum_{i=0}^{n-1} f_0 (T^i x) - 1 \big| < \delta \} 
\to 1$ as $n \to \infty$. 
Hence, 
\[
\lim_{n \to \infty} \mu \big\{ x \in X : \sum_{i=0}^{n-1} f_0 (T^i x) > n ( 1 - \delta ) \big\} = 1 . 
\]
Since $\sum_{i=0}^{n-1} f (T^i x)  \geq. \sum_{i=0}^{n-1} f_0 (T^i x)$ for a.e. $x \in X$,  then 
\[
\lim_{n \to \infty} \mu \big\{ x \in X : \sum_{i=0}^{n-1} f (T^i x) > n ( 1 - \delta ) \big\} = 1 . 
\]
Since Schmidt's condition (\ref{schmidt}) does not hold, there is no measurable solution $g$. 

The final case to prove is where 
$\int_{f > 0} f d\mu = \int_{f < 0} \big( - f \big) d \mu = \infty$. 
It is proved using a construction similar 
to the one used in Theorem \ref{pos-thm}. 
Choose disjoint measurable sets $X_n \subset X$ for $n\in \integers$ such that 
$f$ is bounded on $X_n$, $\int_{X_n} f d \mu = 0$ and 
\[
\mu \Big( \bigcup_{i=0}^{\infty} X_i \Big) = 1 . 
\]
Let $\epsilon_n > 0$ for $n \in \natural$ be such that $\sum_{n=1}^{\infty} \epsilon_n < \infty$. 
Use Lemma \ref{dec2} to construct potentially infinitely many towers such that 
the function $f_A$, which sums the values of $f$ from the bottom to the top of each tower, 
is bounded.  By Proposition \ref{bnd-prop}, $f_A$ is coboundary for an ergodic measure preserving 
transformation $T_A$ which is defined on the bases of the towers. 
Thus, by Lemma \ref{cob-ext-lem}, the extension transformation $T$ has coboundary $f$. 
Also, the explicit transfer function defined by Lemma \ref{cob-ext-lem} is measurable. $\Box$
\end{pf}

\noindent {\bf Acknowledgments}:  We would like to thank Cesar Silva for feedback 
on this paper, and in particular for pointing out applications of coboundaries 
to invariant $\sigma$-finite measures of non-singular transformations. 
Also, we'd like to thank El Houcein El Abdalaoui and 
Matthijs Borst for providing feedback on a previous version. 


\section{Appendix: Coboundary existence for bounded measurable functions}
The main result in this section is Proposition \ref{bnd-prop}.  This result was previously 
proved in \cite{TAJR}, although we include a proof in this appendix, as well as 
Lemma \ref{ucl} which is used in section \ref{COB}. 
\subsection{Balanced Partitions}
Let $A$ be a measurable subset of $X$ and $f:A\to \real$ in $L_1(A,\mu_A)$.
Let $\epsilon > 0$. We say a finite partition $\Pi$ of $A$ is $\epsilon$-balanced and uniform,
if there exists $E\in \Pi$ such that:
\begin{enumerate}
\item $\mu (E) < \epsilon \mu(A)$,
\item $\int_{A\setminus E} f d\mu = \frac{\mu (A\setminus E)}{\mu. (A)} \int_A f d\mu$,
\item $| f(x) - f(y) | < \epsilon$ for $x,y\in a$ and $a\in \Pi\setminus \{ E\}$,
\item $\mu (c) = \mu (d)$ for $c,d \in \Pi\setminus \{ E\}$.
\end{enumerate}
We refer to this type of partition as a PUB($\epsilon$) partition for $f_{|A}$.
The set $E$ is referred to as the exceptional set of the PUB.
\begin{lem}
Suppose $A\subset X$ is measurable and $f:A\to \real$ is integrable and takes
on essentially infinitely many values. Given $\epsilon >0$, 
there exists a PUB($\epsilon$) partition such that $f$ takes on essentially
infinitely many values on both its exceptional set $E$ and its complement
$A\setminus E$.
\end{lem}

\noindent
{\bf Proof:}
Without loss of generality, it is sufficient to prove the lemma where $0 < || f ||_{\infty} < 1$
and $\epsilon < 1$.  Let $N \in \natural$.
Choose $m\in \natural$ such that
\begin{equation}
\frac{2}{m} < \epsilon .
\end{equation}
For $i = 0,1,2,\ldots , 2m - 1$, let
\begin{equation}
A_i = \{ x\in A: -1 + \frac{i}{m} \leq f(x) < -1 + \frac{i+1}{m} \} .
\end{equation}

Let $\alpha = \min{ \{ \mu (A_i) : \mu (A_i) > 0 \} }$.
There exists $i_0$ such that $f$ takes on infinitely many values
on $A_{i_0}$. Let $E_0$ and $E_1$ be disjoint subsets of $A_{i_0}$
with equal measure and such that
\begin{eqnarray}
\frac{1}{\mu (E_0)} \int_{E_0} f d\mu < \frac{1}{\mu (A_{i_0})} \int_{A_{i_0}} f d\mu ,\\
\frac{1}{\mu (E_1)} \int_{E_1} f d\mu > \frac{1}{\mu (A_{i_0})} \int_{A_{i_0}} f d\mu ,
\end{eqnarray}
and $f$ takes on infinitely many values on the set
$A_{i_0} \setminus (E_0 \cup E_1)$ and on the set
$E_0 \cup E_1$.
Let
$$d = \min{\{ | \frac{1}{\mu (E_i)} \int_{E_i} f d\mu - \frac{1}{\mu (A_{i_0})} \int_{A_{i_0}} f d\mu | : i = 0,1 \} } .$$

By simultaneous Diophantine approximation \cite{Cas57},
there exist $q \in \natural$ and $p_i \in \natural$ such that
\begin{eqnarray}
q & > & \max{\{ \frac{2N}{(1 - \epsilon)\mu (A)}, \frac{2\mu (A)}{d\mu (E_1)} \}},
\end{eqnarray}
and for $i=0, 1, \ldots , 2m - 1$,
\begin{eqnarray}
| q \mu (A_i) - p_i | & < & q^{\frac{-1}{2m}}, \\
2m q^{{-1}/{2m}} & < & \epsilon ,\\
2mq^{{-1}/{2m}} & < & d ( \frac{2\alpha}{3} - q^{\frac{-1}{2m}} ) .
\end{eqnarray}

Let $n = q + 1$. Thus,
\begin{eqnarray}
|\mu (A_i) - (\frac{p_i}{n} + \frac{\mu (A_i)}{n})| < n^{-1}q^{{-1}/{2m}}.
\end{eqnarray}
Let $h = \sum\limits_{i=0}^{2m-1} p_i$.
For $i=0, 1, \ldots , 2m - 1$, we can choose subsets
$B_i \subset A_i$ such that
\begin{eqnarray}
\mu (B_i) & = & \mu (A_i) - \frac{p_i}{n} ,\\
\frac{1}{\mu (B_i)} \int_{B_i} f d\mu & = & \frac{1}{\mu (A_i)} \int_{A_i} f d\mu .
\end{eqnarray}
Thus,
\begin{eqnarray}
| \sum\limits_{i=0}^{2m-1} \int_{B_i} f d\mu | & = & | \sum\limits_{i=0}^{2m-1} \frac{\mu (B_i)}{\mu (A_i)} \int_{A_i} f d\mu |
= | \sum\limits_{i=0}^{2m-1} (\frac{\mu (B_i)}{\mu (A_i)} - \frac{1}{n}) \int_{A_i} f d\mu | \\
& \leq & \sum\limits_{i=0}^{2m-1} | \mu (B_i) - \frac{\mu (A_i)}{n} |
= \sum\limits_{i=0}^{2m-1} | \mu (A_i) - \frac{p_i + \mu (A_i)}{n} | \\
& < & 2m n^{-1} q^{{-1} / {2m}} < \frac{d}{n}
(\frac{2\alpha}{3} - q^{\frac{-1}{2m}} ) .
\end{eqnarray}
This implies we can choose $B_{i_0}$ such that
\begin{eqnarray}
\sum\limits_{i=0}^{2m-1} \int_{B_i} f d\mu = 0 .
\end{eqnarray}
Let $E = \bigcup_{i=0}^{2m-1} B_i$ and partition each
set $A_i\setminus B_i$ into $p_i$ subsets of measure ${1}/{n}$
to form $\Pi$.
Therefore, $\mu(E) < \epsilon$ and our lemma is proven. $\Box$

\subsection{Balanced Uniform Towers}
Let $A$ be a measurable subset of $X$ and $f:A\to \real$ a bounded, mean-zero function.
Given finite measurable partition $Q$, 
$h\in \natural$ and $\epsilon > 0$, an $\epsilon$-balanced and uniform tower for $f$ is a set
of disjoint measurable sets $I_i\subset A$ for $i=1,2,\ldots ,h$ and an invertible measure preserving
map $T: I_i \to I_{i+1}$ for $i=1,2,\ldots ,h-1$, such that:
\begin{eqnarray}
\mu (\bigcup_{i=1}^h I_i) & > & \mu (A) - \epsilon , \\
| f( x ) - f( y ) | & < & \epsilon \ \ \mbox{ for } x, y \in I_i, 1 \leq i < h ,\label{nbig0} \\
|\sum\limits_{i=0}^{k} f(T^i x) | & < & \| f \|_{\infty} + \epsilon \ \ \mbox{ for } x\in I_1, k < h ,\label{nbig1}\\
| \sum\limits_{i=0}^{h-1} f(T^i x) | & < & \epsilon \ \ \mbox{ for } x\in I_1 , \label{nbig2} \\ 
\mbox{for each}\ q\in Q, &\exists& \mathcal{I} \subset \{ 1,\ldots ,h \} \ \mbox{such that}\ \mu 
\big( q \triangle ( \bigcup_{i \in \mathcal{I}} I_i ) \big) < \epsilon. 
\end{eqnarray}
We refer to this type of tower as a $TUB(\epsilon, h, Q)$ tower for $f_{|A}$.

\begin{lem}
\label{ucl}
Let $( X, \B , \mu )$ be a standard probability space and $A$ a measurable subset of $X$.
Suppose $f: A \to \real$, $f \in L^{\infty}_0$, takes on essentially infinitely many values. 
Given $N \in \natural$, $\epsilon > 0$ and finite measurable partition $Q$, there exists $h > N$
such that $f$ has a $TUB(\epsilon, h, Q)$ tower.
\end{lem}

\noindent
{\bf Proof:}
From the construction of PUB(${\epsilon} / {3}$) in the previous lemma,
partition $A_i \setminus B_i$ into a disjoint union of sets $A_i(j)$
for $j=1,2,\ldots ,p_i$, such that
\begin{eqnarray}
\mu (A_i(j)) & = & \frac{1}{n} .
\end{eqnarray}

\subsubsection{Greedy Stacking}
Now we give an inductive procedure for stacking the sets $A_i(j)$.
Choose arbitrary $A_i(j)$ and label the set $I_1$.
Given $I_1, I_2, \ldots , I_{k-1}$, let
\begin{equation}
\sigma_{k-1} = \sum\limits_{i=1}^{k-1} \int_{I_i} f d\mu .
\end{equation}
If $k = h$, then we are done.
If $\sigma_{k-1} \leq 0$,
choose
$$
I_k=A_i(j) \not\subset \bigcup_{i=1}^{k-1} \{ I_i \}
$$
such that
$\int_{I_k} f d\mu \geq 0$.
This is possible,
since $k < h $ and $\sigma_h = \sum\limits_{i=0}^{2m-1}\int_{A_i\setminus B_i} f d\mu = 0$.
Otherwise, if $\sigma_k > 0$,
then by the construction of $A_i(j)$, there
exists $I_k \not\subset \bigcup_{i=1}^{k-1} \{ I_i \}$ such that
$\int_{I_k} f d\mu < 0$.
This procedure produces a sequence of sets $I_i$ for $i=1, 2, \ldots , h$
with the property:
\begin{eqnarray}
\sum\limits_{i=1}^{h} \int_{I_i} f d\mu & = & \sum\limits_{i=0}^{2m-1} \int_{A_i\setminus B_i} f d\mu \\
& = & \sum\limits_{i=0}^{2m-1} \int_{B_i} f d\mu = 0 .
\end{eqnarray}

\subsubsection{Level Refinement}
Our transformation $T$ will map $I_i$ onto $I_{i+1}$ for $i = 1, 2, \ldots , h-1$.
Let $\Phi$ be the set of measure preserving maps $T$ such that
$I_{i+1} = T(I_i)$ for $i=1,2,\ldots ,h-1$.
Given $T\in \Phi$, disjoint subsets $D_1, D_2$ contained in $I_1$ with equal measure,
and an invertible measure preserving mapping $\psi : D_1 \to D_2$, let
\[
d(T, D_1,D_2, \psi) = \inf_{x\in D_1} ( \sum\limits_{i=0}^{h-1} f(T^ix) - \sum\limits_{i=0}^{h-1} f(T^i(\psi (x))) ) .
\]
Define
\[
d(D_1,D_2) = \sup_{\psi} d(T, D_1,D_2,\psi ) .
\]
and
\[
d(T) = \sup_{D_1} \{ \mu (D_1) :\, \text {there exists}\, D_2 \mbox{ such that } d(D_1,D_2) > \epsilon \} .
\]
Finally, let
\[
d = \inf_{T\in \Phi} d(T) .
\]
We claim that $d=0$. If $d > 0$, then there exists $T\in \Phi$ such that
$| d(T) - d | < {d} / {h}$. This produces $D_1, D_2$ and $\psi$
such that $d(T,D_1,D_2,\psi ) \geq \epsilon$ and $| \mu (D_1) - d | < {d} / {h}$.
Then there exists $0\leq i < h$ such that for $x\in D_1$,
\[
f(T^ix) \leq f(T^i(\psi x)) - \frac{\epsilon}{h} .
\]
Modify the map $T$, by switching $T^i(D_1)$ and $T^i(D_2)$. Thus, there exists
$T_1 \in \Phi$ such that $T_1 (T^{i-1}D_1) = T^i(D_2)$ and
$T_1 (T^{i-1}D_2) = T^i(D_1)$.
If $d(T_1, D_1,D_2,\psi) \geq \epsilon$, modify $T_1$ in a similar manner
to produce $T_2$. After a finite number of steps, we may produce
$T_k$ such that $d(T_k,D_1,D_2,\psi ) < \epsilon$.
By passing to a subset of $D_1$ if necessary, we obtain
$T^{\prime} \in \Phi$ such that $d(T^{\prime}) < d$ which proves that $d=0$ by contradiction.
Therefore, this proves (\ref{nbig2}) of our lemma. Claim (\ref{nbig1})
follows in a similar manner. $\Box$
\bigskip

\subsubsection{Level Refinement (alternative proof)}
Our transformation $\tau$ will map $I_i$ onto $I_{i+1}$ for $i = 1, 2, \ldots , h-1$. 
Choose $k$ such that $k > {3h} / {\epsilon}$.  
Partition the range such that for $\ell =-k, -k + 1, \ldots , -1, 0, 1, \ldots , k - 1$, 
\[
B_{i, \ell} = \{ x \in I_i : \frac{\ell}{k} \leq f(x) < \frac{\ell+1}{k} \} . 
\]
WLOG, assume $X$ is an ordered set (i.e., $[0,1]$). 
Let $\tau_0$ be any invertible measure preserving map such that $\tau_0 : I_i \to I_{i+1}$ 
for $i \in \{ 1, 2, \ldots , h-1 \}$.  
Define the quantized function $f_i ( x ) = { \ell } / { k }$ if $x \in B_{ i, \ell }$. 
Define an invertible measure preserving map 
$\psi_i : I_i \to I_i$ such that $f_i \circ \psi_i$ is non-decreasing for $i \geq 2$. 
Let $\psi_1$ be the identity map on $I_1$.  
Define an invertible measure preserving map 
$\phi_1 : I_1 \to I_1$ such that $f_1 \circ \phi_1$ is non-increasing. 
The map $g_2 = f_1 \circ \phi_1 + f_2 \circ \psi_2 \circ \tau_0$ is a step function on $I_1$. 
Thus, there exists $\phi_2 : I_1 \to I_1$ such that $g_2 \circ \phi_2$ is non-increasing. 
Let $g_3 = g_2 \circ \phi_2 + f_3 \circ \psi_3 \circ \tau_0^2$.  Continue this process until we have defined $g_h$. 
In particular, by induction, $g_{h-1}$ will be a step function on $I_1$. 
Thus, we can define $\phi_{h-1}: I_1 \to I_1$ such that $g_{h-1} \circ \phi_{h-1}$ 
is non-increasing.  Let $g_h = g_{h-1} \circ \phi_{h-1} + f_h \circ \psi_h \circ \tau_0^{h-1}$.  
A formula for $g_h$ is 
\[
g_h = \sum_{\ell =1}^{h} f_{\ell} \psi_{\ell} \tau_0^{\ell -1} \Pi_{j=\ell}^{h-1} \phi_j . 
\]
For $2 \leq \ell \leq h$, define 
\[
\tau_{\ell} = \psi_{\ell} \tau_0^{\ell - 1} \Pi_{j=\ell}^{h-1} \phi_j . 
\]
Let $\tau_1$ be the identity map. 
Each $\tau_{\ell}$ is an invertible measure preserving mapping 
from $I_1 \to I_{\ell}$.  Define the final mapping $\tau$ as 
$\tau ( x ) = \tau_{\ell +1} \circ \tau_{\ell}^{-1} (x)$ for $x \in I_{\ell}$. 
Because of the greedy algorithm of sorting at each stage and re-ordering so that the next 
level has $f_{\ell}$ monotonic in the opposite direction, then 
the quantized functions $f_{\ell}$ do not exhibit much variation 
as points are iterated through the TUB under $\tau$. 
\begin{claim}
For $m \in \natural$, $m < h$ and a.e. $x, y \in I_1$, 
\[
\left| \sum_{\ell =1}^{m} f_{\ell} ( \tau^{\ell -1} x ) - \sum_{\ell =1}^{m} f_{\ell} ( \tau^{\ell -1} y ) \right| < \frac{\epsilon}{3} . 
\]
I.e., there does not exist a pair of sets $D_1, D_2$ of positive measure such that for $x \in D_1$ 
and $y\in D_2$, 
\[
\left| \sum_{\ell =1}^{m} f_{\ell} ( \tau^{\ell -1} x ) - \sum_{\ell =1}^{m} f_{\ell} ( \tau^{\ell -1} y ) \right| \geq  \frac{\epsilon}{3} . 
\]
\end{claim}
\begin{pf}
It is sufficient to prove that for $m \in \natural$, $m < h$ and a.e. $x, y \in I_1$, 
\[
\Big( \sum_{\ell =1}^{m} f_{\ell} ( \tau^{\ell -1} x ) - \sum_{\ell =1}^{m} f_{\ell} ( \tau^{\ell -1} y ) \Big) 
< \frac{\epsilon}{3} . 
\]
By applying the invertible measure preserving isomorphism, 
$\phi_{h-1}^{-1} \phi_{h-2}^{-1} \ldots \phi_{m}^{-1}$, it is sufficient to prove 
for $m \in \natural$, $m < h$ and a.e. $x, y \in I_1$, 
\[
\Big( \sum_{\ell =1}^{m} f_{\ell} \psi_{\ell} \tau_0^{\ell -1} \Pi_{i=\ell}^{m-1} \phi_i (x) - 
\sum_{\ell =1}^{m} f_{\ell} \psi_{\ell} \tau_0^{\ell -1} \Pi_{i=\ell}^{m-1} \phi_i (y) \Big) < \frac{\epsilon}{3} . 
\]
We can prove the claim inductively on $m$.  Clearly, it is true for $m=1$ (by applying the PUB 
condition on $I_1$). 
Suppose it is true for $m < h$. 
Let $x_0$ and $y_0$ be distinct points in $I_1$. 
Let $x_1 = \phi_m(x_0)$ and $y_1 = \phi_m(y_0)$. 
Consider first the case:
\[
0 < \sum_{\ell =1}^{m} f_{\ell} \psi_{\ell} \tau_0^{\ell -1} \Pi_{i=\ell}^{m-1} \phi_i (x_1) - 
\sum_{\ell =1}^{m} f_{\ell} \psi_{\ell} \tau_0^{\ell -1} \Pi_{i=\ell}^{m-1} \phi_i (y_1)  < \frac{\epsilon}{3} . 
\]
By the construction of $\phi_m$, $x_0 < y_0$. 
This is because the following function is non-increasing in $x$:
\[
\sum_{\ell =1}^{m} f_{\ell} \psi_{\ell} \tau_0^{\ell -1} \Pi_{i=\ell}^{m} \phi_i (x) . 
\]
Since the function, 
\[
f_{m+1} \psi_{m+1} \tau_0^{m} (x)
\]
is non-decreasing in $x$, then 
\[
f_{m+1} \psi_{m+1} \tau_0^{m} (x_0) \leq f_{m+1} \psi_{m+1} \tau_0^{m} (y_0) . 
\]
By combining terms, 
\[
\sum_{\ell =1}^{m+1} f_{\ell} \psi_{\ell} \tau_0^{\ell -1} \Pi_{i=\ell}^{m} \phi_i (x_0) - 
\sum_{\ell =1}^{m+1} f_{\ell} \psi_{\ell} \tau_0^{\ell -1} \Pi_{i=\ell}^{m} \phi_i (y_0)  < \frac{\epsilon}{3} . 
\]
The case where 
\[
0 < \sum_{\ell =1}^{m} f_{\ell} \psi_{\ell} \tau_0^{\ell -1} \Pi_{i=\ell}^{m-1} \phi_i (y_1) - 
\sum_{\ell =1}^{m} f_{\ell} \psi_{\ell} \tau_0^{\ell -1} \Pi_{i=\ell}^{m-1} \phi_i (x_1)  < \frac{\epsilon}{3} , 
\]
may be handled in a similar fashion.  This completes the proof of the claim. 
$\Box$ 
\end{pf}

\noindent
Now we complete the proof of the lemma.  
The function $f$ was quantized to $f_{\ell}$ in such a way that 
for $x \in I_1$ and $m \in \{ 1, 2, \ldots , h \}$, 
\begin{eqnarray}
| \sum_{\ell=1}^{m} f ( \tau^{\ell-1} x ) - \sum_{\ell=1}^{m} f_{\ell} ( \tau^{\ell-1} x ) | 
&\leq& \sum_{\ell=1}^{m} | f ( \tau^{\ell-1} x ) - f_{\ell} ( \tau^{\ell-1} x ) | \\ 
&<& h \Big( \frac{\epsilon}{3h} \Big) = \frac{\epsilon}{3} . 
\end{eqnarray}
Hence, 
\begin{eqnarray}
| \sum_{\ell=1}^{m} f ( \tau^{\ell-1} x ) - \sum_{\ell=1}^{m} f ( \tau^{\ell-1} y ) |  & \leq & 
| \sum_{\ell=1}^{m} f ( \tau^{\ell-1} x ) - \sum_{\ell=1}^{m} f_{\ell} ( \tau^{\ell-1} x ) | \\ 
& + & | \sum_{\ell=1}^{m} f_{\ell} ( \tau^{\ell-1} x ) - \sum_{\ell=1}^{m} f_{\ell} ( \tau^{\ell-1} y ) | \\
& + & | \sum_{\ell=1}^{m} f_{\ell} ( \tau^{\ell-1} y ) - \sum_{\ell=1}^{m} f_{\ell} ( \tau^{\ell-1} y ) | \\
& < & \frac{\epsilon}{3} + \frac{\epsilon}{3} + \frac{\epsilon}{3} = \epsilon , 
\end{eqnarray}
Therefore, this proves (\ref{nbig2}) of our lemma. Claim (\ref{nbig1})
follows in a similar manner. $\Box$
\bigskip

The following proposition was previously proved in \cite{TAJR}.  
For completeness, we provide a similar proof here. 
\begin{prop}
\label{bnd-prop}
Suppose $f: [0,1] \to \real$ is measurable, mean-zero and bounded. 
There exists an ergodic measure preserving transformation $T$ and bounded function $g$ 
such that $f = g - g\circ T$ a.e.  Moreover, the transformation $T$ and 
transfer function $g$ may be constructed 
such that for any $\delta > 0$, $|| g ||_{\infty} < || f ||_{\infty} + \delta$. 
\end{prop}

\noindent
\begin{pf}
If $f = \sum_{i=1}^{m} a_i I_{A_i}$ is a finite step function, a solution is given in \cite{TAJR}.  
The transfer function $g$ is bounded, since, by \cite{LS}, 
\[
g(x) = \lim_{n\to \infty} \frac{1}{n} \sum_{k=1}^{n} \sum_{i=0}^{k-1} f ( T^i x ) \leq \sum_{i=1}^{m} | a_i | . 
\]
Otherwise, $f$ takes on essentially infinitely many values. 
Let $\delta_i > 0$ be such that $\sum_{i=1}^{\infty} \delta_i < \infty$ 
and $Q_i$ for $i\in \natural$ a refining sequence of partitions which 
generate the sigma algebra $\mathcal{B}$. 
Let $\epsilon_1 = \delta_1$. 
Use Lemma \ref{ucl} to construct a $TUB(\epsilon_1, h_1, Q_1)$ tower with decomposition 
into sets $A_1$, $B_1$ and measure preserving map $T_1$. 
Also, assume $A_1$ is made of levels $I_{1,i}$ for $1\leq i \leq h_1$. 
Let $S_1 = T_1$. 
Define $f_1: B_1 \cup I_{1,1} \to \real$ by 
\begin{eqnarray*} 
f_1(x)= 
\left\{\begin{array}{ll}
\sum_{i=0}^{h_1-1} f(S_1^i x) , & \mbox{if $x\in I_{1,1}$}, \\ 
f(x) , & \mbox{if $x\in B_1$} . 
\end{array}
\right.
\end{eqnarray*}
Let $\epsilon_2 = \delta_2 \mu ( I_{1,1} )$. 
Since $\int_{B_1 \cup I_{1,1}} f_1 d\mu = 0$, then we can apply Lemma \ref{ucl} 
to $f_1$ to obtain a $TUB(\epsilon_2, h_2, Q_2)$ tower and 
decompose $B_1 \cup I_{1,1}$ into $A_2 = \cup_{i=1}^{h_2} I_{2,i}$ and $B_2$ 
such that there exists measure preserving $T_2: I_{2,i} \to I_{2,i+1}$ 
for $i = 1, \ldots , h_2 - 1$. 
Define $S_2$ as 
\begin{eqnarray*} 
S_2(x)= 
\left\{\begin{array}{ll}
S_1 (x) , & \mbox{if $x\in S_1^i I_{2,j} \subset S_1^i I_{1,1}$, for $0\leq i \leq h_1-1$ and $1 \leq j \leq h_2$}, \\ 
T_2 ( S_1^{ 1 - h_1 } (x) ) , & \mbox{if $x\in S_1^{h_1 - 1} I_{2,j} \subset S_1^{h_1 - 1} I_{1,1}$, for $1 \leq j < h_2$}, \\ 
T_2(x) , & \mbox{if $x\in B_1 \cap A_2 \setminus I_{2,h_2}$} . 
\end{array}
\right.
\end{eqnarray*}
Suppose $T_n$ and $S_n$ have been defined. 
Proceed in a similar manner to define $S_{n+1}$.  In particular, for a.e. $y \in X\setminus B_n$, 
there exists a unique $x \in I_{n,1}$ and $j_y \geq 0$ such that 
$y = S_n^{j_y} x$.  
For a.e. $x \in I_{n,1}$, 
there exists a minimum $k_{n,x} \geq 0$ such that $S_n^{k_{n,x}} x \in I_{n,h_n}$. 
Define $f_n : B_n \cup I_{n,1} \to \real$ such that 
\begin{eqnarray*} 
f_n(x)= 
\left\{\begin{array}{ll}
\sum_{i=0}^{k_{n,x}} f(S_n^i x) , & \mbox{if $x\in I_{n,1}$}, \\ 
f(x) , & \mbox{if $x\in B_n$} . 
\end{array}
\right.
\end{eqnarray*}
Let $\epsilon_{n+1} = \delta_{n+1} \mu ( I_{n,1} )$. 
Since $\int_{B_n \cup I_{n,1}} f_n d\mu = 0$, then we can apply Lemma \ref{ucl} 
to $f_n$ to obtain a $TUB(\epsilon_{n+1}, h_{n+1}, Q_{n+1})$ tower and 
decompose $B_n \cup I_{n,1}$ into $A_{n+1} = \cup_{i=1}^{h_{n+1}} I_{{n+1},i}$ and $B_{n+1}$ 
such that there exists measure preserving $T_{n+1}: I_{n+1,i} \to I_{n+1,i+1}$ 
for $i = 1, \ldots , h_{n+1} - 1$. 
Define $S_{n+1}$ as 
\begin{eqnarray*} 
S_{n+1}(x)= 
\left\{\begin{array}{ll}
S_n (x) , & \mbox{if $x\in S_n^i I_{n+1,j} \subset S_n^i I_{n,1}$, 
for $0\leq i \leq h_n-1$ and $1 \leq j \leq h_{n+1}$}, \\ 
T_{n+1} ( S_n^{ 1 - h_n } (x) ) , & \mbox{if $x\in S_n^{h_n - 1} I_{n+1,j} \subset S_n^{h_n - 1} I_{n,1}$, for $1 \leq j < h_{n+1}$}, \\ 
T_{n+1}(x) , & \mbox{if $x\in B_n \cap A_{n+1} \setminus I_{n+1,h_{n+1}}$} . 
\end{array}
\right.
\end{eqnarray*}
Note that $S_{n+1} (x) = S_n (x)$ except for $x$ in a set of measure less than or equal to 
\[
\beta_n = \mu \Big( \{ S_n^{k_{n,x}} ( \omega ) : \omega \in I_{n,1} \} \cup 
\{ S_{n+1}^{k_{n+1,x}} ( \omega ) : \omega \in I_{n+1,1} \} 
\cup B_n \cup B_{n+1} \Big) . 
\]
Since $\sum_{n=1}^{\infty} \beta_n < \infty$, 
then $S(x) = \lim_{n\to \infty} S_n (x)$ exists a.e.  
In particular, $S_n(x)$ is eventually constant for a.e. $x \in X$. 
Thus, since each $S_n$ is invertible and measure preserving, 
then $S$ is invertible and measure preserving. 
By a careful choice of $Q_i$, $S$ will be ergodic. $\Box$
\end{pf}


\bigskip

{\small
\parbox[t]{5in}
{T. Adams\\
U.S. Government\\
E-mail: terry@ieee.org\\
}
\smallskip

{\small
\parbox[t]{5in}
{J. Rosenblatt \\
Department of Mathematical Sciences\\
Indiana University-Purdue University Indianapolis\\
Indianapolis, IN 46202, USA\\
E-mail: joserose@iupui.edu}
}
\bigskip


\begin{thebibliography}{99}

\bibitem{Aar97} J. Aaronson, {\em An Introduction to Infinite Ergodic Theory}, Mathematical Surveys and Monographs 50, AMS (1997).

\bibitem{TAJR} T. Adams, J. Rosenblatt, {\em Joint coboundaries}, Dynamical Systems, Ergodic Theory, and Probability: in Memory of Kolya Chernov, Contemporary Mathematics, 698 (2017).

%
%
\bibitem{Anosov73} D.V. Anosov, {\em On an additive functional homology equation connected with an ergodic rotation of the cricle}, Math. USSR Izvestija 7, no. 6 (1973).

%
%
%
%
%
%
%
%
%
%
%

\bibitem{Brow58} F.E. Browder, {\em On the Iteration of Transformations in Noncompact Minimal Dynamical Systems}, Proceedings of the American Mathematical Society, 9, no. 5  (1958)  773-780. 


\bibitem{Cas57}
J.W.S. Cassels, {\em An Introduction to Diophantine Approximation},
Cambridge University Press, London, 1957.

\bibitem{DS2012} A.I. Danilenko and C.E. Silva, {\em Ergodic Theory: Non-singular Transformations},  Mathematics of Complexity and Dynamical Systems. Springer, New York, NY (2012).


\bibitem{DL2001}   Y. Derriennic and M. Lin,
{\em Fractional Poisson equations and ergodic theorems for fractional coboundaries},
Israel Journal of Mathematics 123 (2001) 93-130.

\bibitem{Dot69} W.G. Dotson Jr.,
{\em An application of ergodic theory to the solution of linear functional equations in Banach spaces},
Bull. Amer. Math. Soc. 75 (1969) 347-352.

\bibitem{Dot70} W.G. Dotson Jr.,
{\em On the solution of linear functional equations by averaging iteration},
Proc. Amer. Math. Soc. 25 (1970) 504-506.

\bibitem{Dot71} W.G. Dotson Jr.,
{\em Mean ergodic theorem and iterative solution of linear functional equations},
J. Math. Anal. Appl. 34 (1971) 141-150.

\bibitem{Gro76} C.W. Groftsch, 
{\em Ergodic theory and iterative solution of linear equations},
Applicable Anal. 5 (1976) 313-321.


\bibitem{G2018} D. Giraudo, {\em Invariance principle via orthomartingale approximation}, Stochastics 
and Dynamics 18, no. 6 (2018). 

\bibitem{G69} M.I. Gordin, {\em The central limit theorem for stationary processes}, Sov. Math. Dokl. 10
(1969), 174-176. 

\bibitem{G2009} M.I. Gordin, {\em Martingale-coboundary representation for a class of random fields}, 
Zap. Nauchn. Sem. POMI 364 (2009), 88-108.

\bibitem{GH55} W.H. Gottschalk and G.A. Hedlund, {\em Topological dynamics}, Colloquium publications (American Mathematical Society) 36 (1955).

\bibitem{GO2005} S. Gou\"{e}zel, {\em Regularity of coboundaries for nonuniformly expanding Markov maps}, Proc. of. the AMS 134, no. 2 (2005) 391-401.



\bibitem{Hal76} G. Hal\'{a}sz, {\em Remarks on the remainder in Birkhoff's ergodic theorem}, 
Acta Mathematica Academiae Scientiarum Hungaricae Tomus, 28 (3-4) (1976) 389-395.

\bibitem{KS97} A. Katok and R.J. Spatzier, {\em Differential rigidity of Anosov actions of higher rank abelian groups and algebraic lattice actions}, Tr.  Mat.  Inst.  Steklova, 216  (1997) 292-319.

\bibitem{KS94} A. Katok and R.J. Spatzier, {\em First cohomology of Anosov actions of higher rank abelian groups and applications to rigidity}, Inst. Hautes \'{E}tudes Sci. Publ. Math., 79 (1994) 131-156.

\bibitem{KKM2017} A. Korepanov, Z. Kosloff and I. Melbourne, {\em Martingale-coboundary decomposition for families of dynamical systems}, 
Annales de l'Institut Henri Poincar\'{e} C, Analyse non lin\'{e}aire 35, no. 4 (2018) 859-885.

\bibitem{Ko2002}  I. Kornfeld, {\em Some old and new Rokhlin towers}, Contemporary Mathematics 356, 
Chapel Hill Ergodic Theory Workshops (2004).

\bibitem{KL}  I. Kornfeld and V. Losert, {\em Coboundaries and measure-preserving actions of nilpotent and solvable groups}, Ergodic Theory and Dynamical
Systems, 24 no. 3 (2004) 873-890.

\bibitem{K} I. Kornfeld, {\em Coboundaries for commuting transformations}, Illinois J. Math 43 (1999) 528-539.

\bibitem{Kw} S. Kwapien, {\em Linear Functionals Invariant Under Measure Preserving Transformations}, Math. Nachr. 43 (1984) 175-179.


\bibitem{LS} M. Lin and R. Sine, {\em Ergodic theory and the functional equation $(I- T)x = y$}, J. Operator Theory 10 (1983) no. 1, 153-166.

\bibitem{Liv72} A.N. Liv\u{s}ic, {\em Cohomology of dynamical systems}, Math. USSR Izvestija 6, no. 6 (1972).



\bibitem{Neumann1877} C. Neumann, {\em Untersuchungen \H{u}ber das logarithmische und Newtonsche potential}, Leipzig: Teubner 1877.

%

\bibitem{Pet83} K. Petersen, K. {\em Ergodic Theory} Cambridge University Press, Cambridge, MA. 
(1983) MR0833286 (87i:28002)

%
%
%

\bibitem{Schmidt77} K. Schmidt, {\em Cocycles of ergodic transformation group}, MacMillen Lectures in Math,
Vol. 1. Macmillan Company of India, Ltd., Delhi, 1977.


\bibitem{Vin2017} K. Vinhage, {\em Cocycle rigidity of partially hyperbolic abelian actions with almost rank-one factors}, Ergodic Theory \& Dynamical Systems 1-11. doi:10.1017/etds.2017.119

\bibitem{Vol93} D. Volny, {\em Approximating martingales and the central limit theorem for strictly stationary processes}, Stochastic Processes and their Applications 44 (1993) 41-74. 

\bibitem{VW2004} D. Volny and B. Weiss, {\em Coboundaries in $L_0^{\infty}$}, Ann. I. H. Poincar\'{e}, 
PR 40 (2004) 771-778. 


\end{thebibliography}
\end{document}